\documentclass[a4paper,11pt]{article}
\usepackage[utf8]{inputenc}
\usepackage[T1]{fontenc}
\usepackage{lmodern}
\usepackage[english]{babel}
\usepackage[margin=27mm]{geometry}
\usepackage{amsmath,amssymb,amsthm,mathtools}
\usepackage{mathrsfs}
\usepackage{tikz}
\usetikzlibrary{arrows}
\usepackage{pgfplots}
\pgfplotsset{compat=1.4}
\usepackage{booktabs}
\usepackage{enumitem}
\usepackage{caption}
\usepackage[pdfborder={0 0 0}]{hyperref}

\newtheoremstyle{rigidplain}{\topsep}{\topsep}{\itshape}{}{\normalfont\scshape}{.}{0.5em}{}
\newtheoremstyle{rigidremark}{\topsep}{\topsep}{\normalfont}{}{\normalfont\scshape}{.}{0.5em}{}
\theoremstyle{rigidplain}
\newtheorem{theorem}{Theorem}[section]
\newtheorem{corollary}[theorem]{Corollary}
\newtheorem{proposition}[theorem]{Proposition}
\newtheorem{problem}{Problem}
\theoremstyle{rigidremark}
\newtheorem{remark}[theorem]{Remark}
\newtheorem{example}[theorem]{Example}
\numberwithin{equation}{section}

\newcommand{\Z}{\mathbb Z}
\newcommand{\R}{\mathbb R}
\newcommand{\E}{\mathbb E}
\newcommand{\PP}{\mathbb P}
\newcommand{\eps}{\varepsilon}
\newcommand{\ka}[1]{\mathcal{#1}}
\newcommand{\cA}{\mathcal A}
\newcommand{\cD}{\mathcal D}
\newcommand{\cU}{\mathcal U}
\newcommand{\cH}{\mathcal H}
\newcommand{\cP}{\mathcal P}
\newcommand{\two}{\mathbf 2}
\DeclareMathOperator{\supp}{supp}
\DeclareMathOperator{\sgn}{sgn}
\DeclareMathOperator{\Var}{Var}

\newcommand{\windowpanel}[5]{%
\begin{tikzpicture}[scale=0.4]
  \draw[gray!35] (0,0) grid (4,4);
  \foreach \x/\y in {#1}{%
    \pgfmathtruncatemacro{\ins}{(\x>(#2-1))*(\x<(#2+3))*(\y>(#3-1))*(\y<(#3+3))}%
    \ifnum\ins=1 \fill[black] (\x,\y) rectangle ++(1,1);%
    \else \fill[gray!40] (\x,\y) rectangle ++(1,1);\fi}
  \draw[#5, very thick] (#2,#3) rectangle ++(3,3);
  \node[below, font=\small] at (2,-0.05) {$\eps=(#2,#3)$:\ \textbf{#4}};
\end{tikzpicture}}

\newcommand{\localpanel}[2]{%
\begin{tikzpicture}[scale=0.78]
  \draw[gray!35] (0,0) grid (4,4);
  \foreach \x/\y/\c in {#1}{%
    \fill[\c!30] (\x,\y) rectangle ++(1,1);
    \draw[\c!80!black, thick] (\x,\y) rectangle ++(1,1);
    \foreach \dx in {0,1}{\foreach \dy in {0,1}{%
      \pgfmathtruncatemacro{\px}{\x+\dx}%
      \pgfmathtruncatemacro{\py}{\y+\dy}%
      \pgfmathtruncatemacro{\inn}{(\px>0)*(\px<4)*(\py>0)*(\py<4)}%
      \ifnum\inn=1 \fill[\c!80!black] (\px,\py) circle (4.5pt);\fi}}}
  \draw[dashed, gray] (0.65,0.65) rectangle (3.35,3.35);
  \node[below, font=\small] at (2,-0.1) {#2};
\end{tikzpicture}}

\title{Local fluctuations, majorization, and nontransitivity in cube tilings of $\Z^d$}
\author{Andrzej P. Kisielewicz\\[5pt]
\small Instytut Matematyki, Uniwersytet Zielonog\'orski\\
\small ul. Podg\'orna 50, 65-246 Zielona G\'ora, Poland\\
\small Corresponding author: \texttt{A.Kisielewicz@wmie.uz.zgora.pl}}
\hypersetup{pdfauthor={Andrzej P. Kisielewicz}}
\date{}

\begin{document}
\maketitle

\begin{abstract}
Let $\ka P$ be a partition of $\Z^d$ into translates of the discrete unit cube $\{0,1\}^d$. For a uniformly chosen $\eps\in\{0,1\}^d$, we count the cubes of $\ka P$ entirely contained in the observation window $[0,3]^d+\eps$. The resulting score has mean $(3/2)^d$ for every partition, while its distribution reflects the geometry of the tiling. Comparing two such scores, with the shifts chosen independently, defines a game in which the larger score wins. In dimension $d=4$, we construct a geometric realization of the Steinhaus--Trybu{\l}a paradox: three periodic partitions form a cycle in which each beats the next with probability strictly greater than $1/2$. We describe local configurations by finite cube tiling codes and prove that every resulting list of counts is majorized by that of the simple code $C_d=\{0,2\}^d$. Nevertheless, every non-simple code differing from $C_d$ in only two coordinate directions beats it, with a sharp lower bound on the winning margin. We also show that the ratio of wins to losses of $C_d$ has limit superior at most $1$ as the dimension tends to infinity, uniformly over opposing codes. For $4\Z^d$-periodic tilings, changing the window size can repeatedly reverse the winner, but the winning margin eventually stabilizes separately for side lengths $4m+1$ and $4m+3$.
\end{abstract}

\noindent\textbf{Keywords.} Cube tilings, majorization, local number fluctuations, Steinhaus--Trybu{\l}a paradox.

\section{Introduction}\label{sec:intro}

\begin{figure}[ht]
\centering
\begin{tikzpicture}[scale=0.42]
  \begin{scope}
    \begin{scope}
      \clip (-2.5,-2.5) rectangle (6.5,6.5);
      \draw[gray!30] (-3,-3) grid (7,7);
      \foreach \m in {-2,...,3}{\foreach \n in {-2,...,3}{\fill (2*\m,2*\n) rectangle ++(1,1);}}
    \end{scope}
    \draw[red, very thick] (0,0) rectangle (3,3);
    \draw[blue, very thick, dashed] (1,1) rectangle (4,4);
    \fill[red] (0,0) circle (5pt);
    \node at (2,-3.4) {$\ka F^1$};
  \end{scope}
  \begin{scope}[xshift=12cm]
    \begin{scope}
      \clip (-2.5,-2.5) rectangle (6.5,6.5);
      \draw[gray!30] (-3,-3) grid (7,7);
      \foreach \n in {-2,...,3}{\foreach \m in {-3,...,4}{\fill ({2*\m+\n},{2*\n}) rectangle ++(1,1);}}
    \end{scope}
    \draw[red, very thick] (0,0) rectangle (3,3);
    \draw[blue, very thick, dashed] (1,1) rectangle (4,4);
    \fill[red] (0,0) circle (5pt);
    \node at (2,-3.4) {$\ka F^2$};
  \end{scope}
\end{tikzpicture}
\caption{Portions of the fleets $\ka F^1$ and $\ka F^2$ of~\eqref{eq:two-fleets}; the black squares are the members of the fleets and the dot marks the origin. The window $[0,3]^2$ (solid) contains $4$ cubes of $\ka F^1$ and $3$ cubes of $\ka F^2$. The window $[1,4]^2$ (dashed) contains $1$ and $2$ cubes, respectively.}
\label{fig:fleets}
\end{figure}
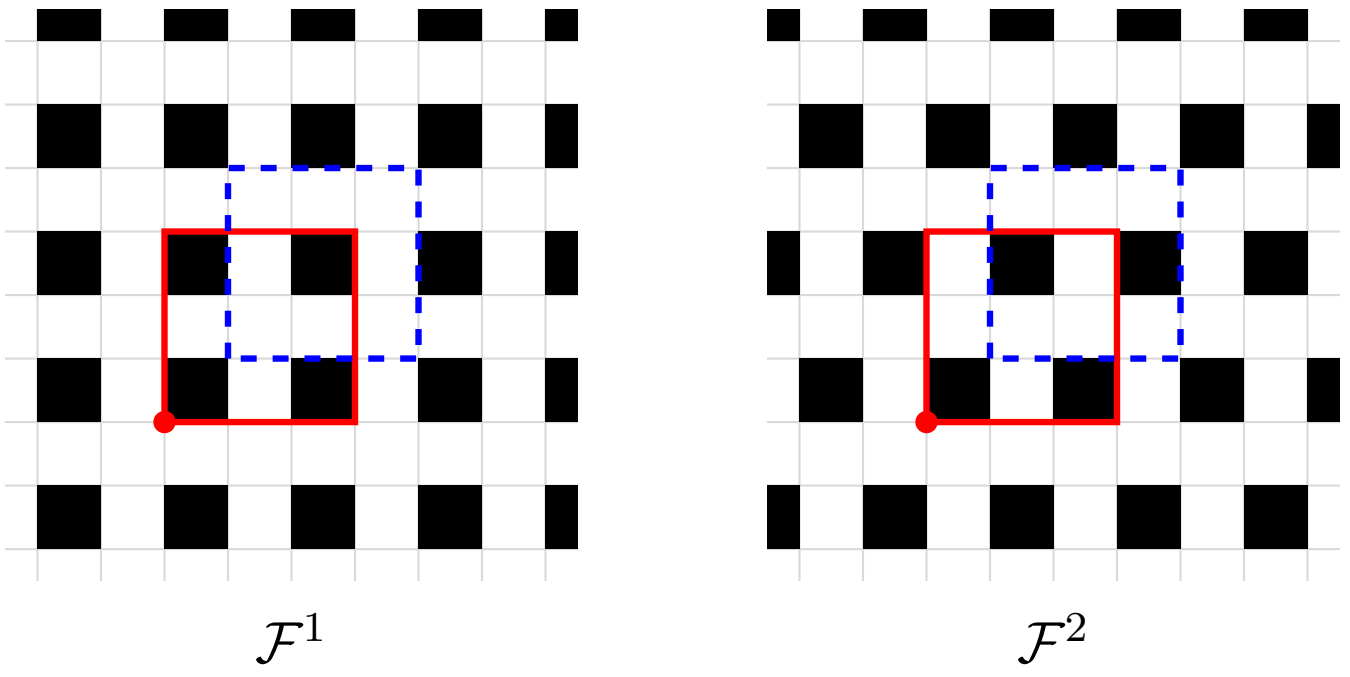

Let $T\subseteq\Z^d$ be a set such that the discrete cubes $\{0,1\}^d+t$, $t\in T$, form a partition of $\Z^d$. The family
\[
\ka F=\{[0,1]^d+t:\ t\in T\}
\]
consists of pairwise disjoint closed unit cubes. We call it a \emph{disjoint fleet}, or simply a \emph{fleet}, and its members \emph{black cubes}. Equivalently, the cubes $[0,2)^d+t$, $t\in T$, tile $\R^d$. Thus fleets are in one-to-one correspondence with tilings of $\R^d$ by translates of $[0,2)^d$ with integer translation vectors. Figure~\ref{fig:fleets} shows portions of the planar fleets
\begin{equation}\label{eq:two-fleets}
\ka F^1=\{[0,1]^2+(2m,2n):m,n\in\Z\}\quad\text{and}\quad 
\ka F^2=\{[0,1]^2+(2m+n,2n):m,n\in\Z\}.
\end{equation}

\paragraph{Observation windows and the game.}
Fix $v\in\Z^d$ and let $E_d=\{0,1\}^d$. For $\eps\in E_d$ the \emph{observation window} is
\[
R(v+\eps)=[0,3]^d+v+\eps .
\]
The \emph{score} of $\ka F$ in this window is the number $N_v(\ka F,\eps)$ of black cubes contained in $R(v+\eps)$. The \emph{list of results} $L_v(\ka F)$ is the list of the $2^d$ scores $N_v(\ka F,\eps)$, $\eps\in E_d$, arranged in nondecreasing order with multiplicities retained. We use this convention for all lists of results.

Two players $A$ and $B$ choose fleets $\ka F_A$ and $\ka F_B$, respectively. Player $A$ chooses a uniform shift $\eps\in E_d$, and player $B$ independently chooses a uniform shift $\delta\in E_d$. The larger score wins, and equal scores give a draw. Let
\[
w_v(\ka F_A,\ka F_B)=\bigl|\{(\eps,\delta)\in E_d^2:N_v(\ka F_A,\eps)>N_v(\ka F_B,\delta)\}\bigr|
\]
be the number of pairs of shifts won by $A$. We say that $\ka F_A$ \emph{beats} $\ka F_B$ (at $v$) if $w_v(\ka F_A,\ka F_B)>w_v(\ka F_B,\ka F_A)$. Players $A$ and $B$ win with probabilities $w_v(\ka F_A,\ka F_B)/4^d$ and $w_v(\ka F_B,\ka F_A)/4^d$, respectively.

\begin{figure}[ht]
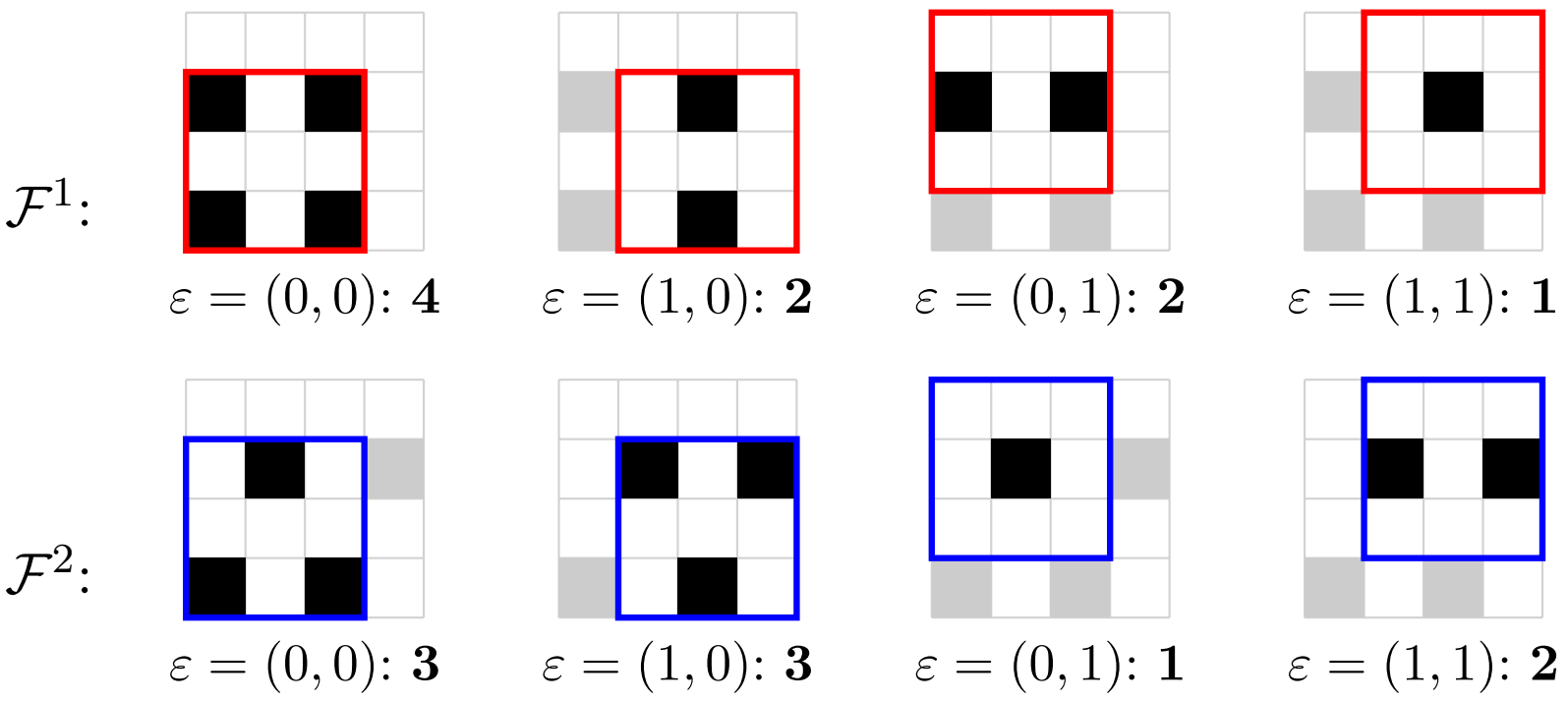

\centering
\begin{tabular}{r@{\quad}cccc}
\raisebox{0.8cm}{$\ka F^1$:} &
\windowpanel{0/0,2/0,0/2,2/2}{0}{0}{4}{red} &
\windowpanel{0/0,2/0,0/2,2/2}{1}{0}{2}{red} &
\windowpanel{0/0,2/0,0/2,2/2}{0}{1}{2}{red} &
\windowpanel{0/0,2/0,0/2,2/2}{1}{1}{1}{red}\\[2pt]
\raisebox{0.8cm}{$\ka F^2$:} &
\windowpanel{0/0,2/0,1/2,3/2}{0}{0}{3}{blue} &
\windowpanel{0/0,2/0,1/2,3/2}{1}{0}{3}{blue} &
\windowpanel{0/0,2/0,1/2,3/2}{0}{1}{1}{blue} &
\windowpanel{0/0,2/0,1/2,3/2}{1}{1}{2}{blue}
\end{tabular}
\caption{The four windows $[0,3]^2+\eps$ inside $[0,4]^2$. Only the black cubes contained in $[0,4]^2$ are drawn; there are always $2^d=4$ of them, and they form the local code of Section~\ref{sec:local}. Cubes inside the window are black, the others grey. The scores are $4,2,2,1$ for $\ka F^1$ and $3,3,1,2$ for $\ka F^2$.}
\label{fig:windows}
\end{figure}

For the fleets~\eqref{eq:two-fleets} and $v=0$ we have (see Figure~\ref{fig:windows})
\[
L_0(\ka F^1)=[1,2,2,4]\quad\text{and}\quad L_0(\ka F^2)=[1,2,3,3].
\]
Thus $w_0(\ka F^1,\ka F^2)=6$ and $w_0(\ka F^2,\ka F^1)=7$, while three of the sixteen pairs of shifts give a draw. Hence $\ka F^2$ beats $\ka F^1$, although both expected scores are equal to $9/4$.

\paragraph{The common sum and the choice of the window.}
For every fleet and every $v$, the $2^d$ scores sum to $3^d$ (Proposition~\ref{prop:count}), so their mean is always $(3/2)^d$. How this sum is distributed among the individual windows depends on the geometry of the fleet. The window side $3$ is the smallest one that detects this geometry. For every fleet and every fixed $v$, exactly one of the $2^d$ windows $[0,1]^d+v+\eps$, where $\eps\in\{0,1\}^d$, contains a black cube, while the windows $[0,2]^d+v$ and $[0,4]^d+v$ always contain exactly $1$ and $2^d$ black cubes, respectively (Proposition~\ref{prop:even-odd-windows}). Windows of sides $1$, $2$ and $4$ therefore give the same list for every fleet. We examine the effect of changing the window size in more detail in Section~\ref{sec:other-windows}.

Deviations of the number of objects in a randomly placed window from its mean are called \emph{local number fluctuations}. Such fluctuations are studied in statistical physics as a way of describing the spatial organization of many-particle systems. Torquato and Stillinger~\cite{TorquatoStillinger2003} relate local density fluctuations and hyperuniformity to the variance of the number of points in a window. Torquato, Kim and Klatt~\cite{TorquatoKimKlatt2021} consider higher moments and the full distribution of this number, and Kim and Torquato~\cite{KimTorquato2017} study the dependence on the shape and orientation of the window. In our model we count black cubes entirely contained in the window, the shifts are taken from a discrete set, and two independent scores are compared.

Another interpretation is to view the black cubes as a fleet of ships and the window as the detection region of a radar. The aim is then to arrange the ships in a formation that minimizes the number detected by the radar. Replacing the rule ``more'' by ``less'' interchanges wins and losses and preserves draws, so the same calculations apply with the winning direction reversed.

Observing fleets through moving windows can also reveal arithmetic structures. In \cite{KeplerCubeDraft}, the scores of a particular lattice fleet, assigned to the corresponding vertices $\eps\in\{0,1\}^d$, form a \emph{Kepler cube} encoding the first $d$ levels of the Kepler tree and their parent--child relations. Face partitions induced by the same lattice encode the first $d$ levels of the Calkin--Wilf tree; allowing $d$ to vary yields the whole tree. These constructions illustrate how local observations of fleets can lead to further geometric and arithmetic results.

\paragraph{Geometrically generated dice.}
The list $L_v(\ka F)$ may be regarded as a die with $2^d$ equally likely faces, one for each shift of the window. All these dice have the same sum of face values $3^d$, so the game is a game of equal-sum dice in the sense of~\cite{FinkelsteinThorp2007,DeSchuymerDeMeyerDeBaets2006}. Nontransitive dice are a classical manifestation of the Steinhaus--Trybu{\l}a paradox~\cite{SteinhausTrybula1959,Trybula1961,Gardner1970}. Here the face values are determined by geometry: they count the unit cubes contained in translated windows. Thus cube tilings produce a constrained family of equal-mean distributions. In Section~\ref{sec:geometric-paradox} we construct three fleets in $\R^4$ whose independent scores satisfy the three cyclic inequalities with probabilities strictly greater than $1/2$. This gives a geometric realization of the paradox without conditioning on the absence of a draw.

\paragraph{Relation to earlier work.}
Dutour Sikiri\'c, Itoh and Poyarkov~\cite[Section~4]{DutourItohPoyarkov2007} study the number of cubes of side $2$ contained in a window of side $4$ for $4\Z^d$-periodic tilings with integer translation vectors. Averaging over all $4^d$ shifts modulo $4$, they show that the regular tiling maximizes the second moment. The counts agree with ours when the same translation vectors are used: in both models, containment requires the translation vector to lie in the same set of three consecutive integers in each coordinate.

Our Theorem~\ref{thm:majorization} strengthens their second-moment bound to a local majorization statement: for every fixed $v$, it compares the entire list of $2^d$ scores, and hence the mean of every convex function of the score. This also applies to nonperiodic fleets through their local codes. For a periodic tiling, taking the square function and averaging over the blocks $2\delta+E_d$, $\delta\in E_d$, recovers the earlier second-moment bound. The further questions addressed here concern comparisons of two independent scores. These include nontransitivity, bounds on winning margins, comparisons of Cartesian powers, and the uniform asymptotic bound of Theorem~\ref{thm:simple-asymptotic}; they are not considered in~\cite{DutourItohPoyarkov2007}.

\paragraph{Main results.}
\begin{enumerate}[label=(\alph*),leftmargin=*]
\item \emph{Local structure} (Section~\ref{sec:codes}). The black cubes of a fleet contained in $[0,4]^d+v$ form, after translation, a cube tiling code $F_v$ over the alphabet $\{0,1,2,3\}$ with complementary pairs $\{0,2\}$ and $\{1,3\}$. We call such a code a \emph{fleet code}. Equivalently, the intersection of the fleet with $\{1,2,3\}^d+v$ is a minimal partition of this box in the sense of~\cite{GrytczukKisielewiczPrzeslawski2004,KisielewiczPrzeslawski2008}. Conversely, every fleet code is the local code of a periodic fleet (Proposition~\ref{prop:local-code}). Each local comparison of fleets, periodic or not, is therefore a comparison of fleet codes.
\item \emph{Majorization} (Section~\ref{sec:majorization}). For every fleet code $F$ in dimension $d$, the list of results of $F$ is majorized by the list of the simple code $C_d=\{0,2\}^d$ (Theorem~\ref{thm:majorization}). Hence $C_d$ maximizes $\E[\varphi(X)]$ for every convex $\varphi$, where $X(j)=f_j$ for $L(F)=[f_1,\ldots,f_{2^d}]$ and $j$ chosen uniformly at random from $[2^d]$, and minimizes the geometric mean of the scores (Corollaries~\ref{cor:convex-order} and~\ref{cor:log-mean}).
\item \emph{Codes beating the simple code} (Section~\ref{sec:two-coordinates}). Every non-simple fleet code that differs from $C_d$ only in two coordinate directions beats $C_d$, and the number of winning pairs exceeds the number of losing pairs by at least $d-1$ (Theorem~\ref{thm:two-coordinate-changes}). The bound is sharp in every dimension $d\ge2$ (Example~\ref{ex:sharpness}). The proof relies on a weighted inequality for down-sets in the Boolean lattice (Proposition~\ref{prop:downsets}), obtained from the local LYM inequality.
\item \emph{Small dimensions} (Section~\ref{sec:low-dimensions}). Exhaustive computations for $d\le4$ show that unbeaten codes exist for $d\le2$ but not for $d=3,4$. The game is nontransitive even along a chain of the majorization order, and it is not invariant under isomorphisms of codes. We also give an explicit geometric realization of the Steinhaus--Trybu{\l}a paradox in dimension $d=4$, with all three cyclic winning probabilities greater than $1/2$ (Proposition~\ref{prop:geometric-paradox}).
\item \emph{Large dimensions} (Section~\ref{sec:large-dimensions}). For codes with distinct mean logarithms of the scores, sufficiently high Cartesian powers are ordered by these means: the power of the code with the larger mean beats that of the code with the smaller mean (Theorem~\ref{thm:lln-products}). In particular, such powers eliminate nontransitive cycles among these codes. The simple code can beat other codes by an arbitrarily large number of pairs (Proposition~\ref{prop:unbounded-simple-advantage}). Nevertheless, the ratio of its wins to its losses has limit superior at most $1$ (Theorem~\ref{thm:simple-asymptotic}).
\item \emph{Other windows} (Section~\ref{sec:other-windows}). We prove counting identities for arbitrary rectangular windows with integer sides. We show that changing the window size can reverse the winner and that the comparison eventually stabilizes separately for side lengths $4m+1$ and $4m+3$. All four types of eventual behavior occur already in dimension $3$, and a four-dimensional example reaches its eventual balances only from the odd side length $53$ onward.
\end{enumerate}
Section~\ref{sec:problems} collects open problems.

\paragraph{Usage of AI.}
 All theorem statements and definitions were formulated by the author. The proofs of Theorems~\ref{thm:two-coordinate-changes} and~\ref{thm:simple-asymptotic} were developed with assistance from GPT Astra~6 Pro. Computations were carried out with assistance from Codex, and the manuscript was prepared and edited with assistance from the AI assistants Astra and Claude. The author takes full responsibility for the results and the entire content of this paper.

\section{Fleet codes}\label{sec:codes}


\subsection{Cube tiling codes and minimal partitions}\label{sec:tiling-codes}

For an integer $r\ge0$, write $[r]=\{1,\ldots,r\}$, with $[0]=\emptyset$. For a binary word $\eps\in\{0,1\}^r$, let $|\eps|=\sum_{i=1}^r\eps_i$.

Let $\cA$ be a finite alphabet with a fixed-point-free involution $a\mapsto a'$. Two words $u,t\in\cA^d$ are \emph{dichotomous} if $u_i=t_i'$ for some $i\in[d]$. A \emph{cube tiling code} is a set $V\subseteq\cA^d$ of $2^d$ pairwise dichotomous words; for example, $\{aa,a'a,ba',b'a'\}$ is a cube tiling code in dimension $d=2$. If $V$ is a cube tiling code, then for every $v\in V$ and every $i\in [d]$ there exists $u\in V\setminus \{v\}$ such that $v_i=u'_i$. For background see~\cite{KisielewiczTowards,KisielewiczPrzeslawski2008,LagariasShor1994}. Two cube tiling codes $V,W\subseteq\cA^d$ are \emph{isomorphic} if there are a permutation $\sigma$ of $[d]$ and bijections $\varphi_i:\cA\to\cA$ with $\varphi_i(a')=\varphi_i(a)'$ for all $a$, such that
\[
W=\{\varphi_1(v_{\sigma(1)})\cdots\varphi_d(v_{\sigma(d)}):v_1\cdots v_d\in V\}.
\]

Let $X=X_1\times\cdots\times X_d$, where $X_i$ are finite sets with $|X_i|\ge2$. A nonempty box $B=B_1\times\cdots\times B_d\subseteq X$ is \emph{proper} if $B_i\neq X_i$ for every $i$. A partition of $X$ into $2^d$ proper boxes is a \emph{minimal partition}~\cite{GrytczukKisielewiczPrzeslawski2004,KisielewiczPrzeslawski2008}. By~\cite[Theorem~2]{GrytczukKisielewiczPrzeslawski2004}, a partition of $X$ into proper boxes is minimal if and only if any two of its boxes $B,D$ are dichotomous, that is, $B_i=X_i\setminus D_i$ for some $i$. If $f_i:\cA\to2^{X_i}\setminus\{\emptyset,X_i\}$ satisfy $f_i(a')=X_i\setminus f_i(a)$, then for every cube tiling code $V$ the boxes $f(u)=f_1(u_1)\times\cdots\times f_d(u_d)$, $u\in V$, form a minimal partition of $X$~\cite{KisielewiczPrzeslawski2008}.

From now on we use the alphabet
\[
\cA=\{0,1,2,3\},\qquad 0'=2,\quad 2'=0,\quad 1'=3\quad\text{and}\quad 3'=1,
\]
that is, $a'\equiv a+2\pmod 4$. (The representation by words over $\{0,1,2,3\}$ appears in the work of Corr\'adi and Szab\'o~\cite{CorradiSzabo1990}, while Lagarias and Shor~\cite{LagariasShor1994} explicitly use the terminology of cube-tiling codes.) A cube tiling code over this alphabet is called a \emph{fleet code}, and $\mathscr F_d$ denotes the set of all fleet codes in dimension $d$ (with fixed coordinate and letter labels). A fleet code $F$ is \emph{simple} if $F=A_1\times\cdots\times A_d$ with $A_i\in\bigl\{\{0,2\},\{1,3\}\bigr\}$ for every $i$. We write
\[
C_d=\{0,2\}^d
\]
and $\two_n=2\cdots2$ for the word consisting of $n$ letters $2$; $C_0=\{\two_0\}$ consists of the empty word.

Two words $x,y$ of a fleet code form an \emph{$i$-twin pair} if they differ only in the $i$th coordinate, where they carry complementary letters. The next two propositions recall standard facts and are stated without proof. The first is the flip described in~\cite[Section~2]{DutourItohPoyarkov2007}, expressed in terms of twin pairs (see also~\cite{KisielewiczGluingCutting}).

\begin{proposition}\label{prop:flip}
Let $F$ be a fleet code and let $\{x,y\}\subseteq F$ be an $i$-twin pair with $\{x_i,y_i\}=\{a,a'\}$. Let $\{b,b'\}$ be the other complementary pair, and let $\tilde x,\tilde y$ be obtained from $x,y$ by replacing the $i$th letters by $b$ and $b'$, respectively. Then $(F\setminus\{x,y\})\cup\{\tilde x,\tilde y\}$ is a fleet code.
\end{proposition}

Cube tiling codes in dimensions $d\le4$ have been classified up to isomorphism over arbitrary alphabets with complementation; see~\cite[Section~3.3]{KisielewiczCodes} and the references therein. For our four-letter alphabet, the classification for $d\le4$ was first given in~\cite[Section~2]{DutourItohPoyarkov2007}. We recall the planar description.

\begin{proposition}\label{prop:planar}
Let $F\subset\{0,1,2,3\}^2$ be a fleet code. Then, up to isomorphism,
\[
F=\{00,02,20,22\}\quad\text{or}\quad F=\{00,02,21,23\}.
\]
\end{proposition}

\subsection{The local code of a fleet}\label{sec:local}

Let $\ka F=\{[0,1]^d+t:t\in T\}$ be a fleet and $v\in\Z^d$. The black cubes contained in $[0,4]^d+v$ are those with $t-v\in\{0,1,2,3\}^d$. Let
\[
F_v=(T-v)\cap\{0,1,2,3\}^d .
\]
A discrete cube $\{0,1\}^d+t$ intersects $X_v=\{1,2,3\}^d+v$ if and only if $t-v\in\{0,1,2,3\}^d$. After translating by $-v$, the coordinate sections of these intersections are
\begin{equation}\label{eq:local-realization}
I_0=\{1\},\qquad I_1=\{1,2\},\qquad I_2=\{2,3\}\quad\text{and}\quad I_3=\{3\}.
\end{equation}
Note that $I_0=X_1\setminus I_2$ and $I_1=X_1\setminus I_3$, where $X_1=\{1,2,3\}$.

\begin{proposition}\label{prop:local-code}\leavevmode
\begin{enumerate}[label=\textup{(\alph*)}]
\item For every fleet $\ka F$ and every $v\in\Z^d$, the family $\cP_v=\{I_{u_1}\times\cdots\times I_{u_d}:u\in F_v\}$ is a minimal partition of $\{1,2,3\}^d$, and $F_v$ is a fleet code.
\item Conversely, for every fleet code $H$ the family
\[
\ka F_H=\{[0,1]^d+h+4z:\ h\in H,\ z\in\Z^d\}
\]
is a fleet, and its local code at the origin is $H$.
\end{enumerate}
\end{proposition}

\begin{proof}
(a) The family of discrete cubes $\{0,1\}^d+T$ is a partition of $\Z^d$, so $\cP_v$ is, after translation, a partition of $\{1,2,3\}^d$ into boxes, which are proper because each $I_a$, where $a\in\{0,1,2,3\}$, is a nonempty proper subset of $\{1,2,3\}$. Each $I_a$ contains exactly one element of $\{1,3\}$, so each box of $\cP_v$ contains exactly one point of $\{1,3\}^d$. Hence $|F_v|=|\cP_v|=2^d$ and $\cP_v$ is a minimal partition. By~\cite[Theorem~2]{GrytczukKisielewiczPrzeslawski2004}, for distinct $u,t\in F_v$ there is $i$ with $I_{u_i}=\{1,2,3\}\setminus I_{t_i}$. The only complementary pairs among the sets~\eqref{eq:local-realization} are $I_0,I_2$ and $I_1,I_3$, so $u_i=t_i'$. Thus $F_v$ is a fleet code.

(b) If $h,h'\in H$ are distinct, then $h'_i\equiv h_i+2\pmod 4$ for some $i$, so the sets $\{h_i,h_i+1\}$ and $\{h'_i,h'_i+1\}$ are disjoint modulo $4$. Hence the $2^d$ sets $\{0,1\}^d+h$, $h\in H$, are pairwise disjoint modulo $4$; each has $2^d$ elements, so together they cover $(\Z/4\Z)^d$. Therefore the family of discrete cubes $\{0,1\}^d+h+4z$, where $h\in H$ and $z\in\Z^d$, is a partition of $\Z^d$. Finally, $h+4z\in\{0,1,2,3\}^d$ if and only if $z=0$.
\end{proof}

The scores at $v$ depend only on $F_v$, since every black cube contained in a window $R(v+\eps)$ lies in $[0,4]^d+v$. By Proposition~\ref{prop:local-code}, every pair of fleet codes arises in this way. From now on we therefore take $v=0$ and work with fleet codes. Our results also apply to nonperiodic fleets, since each local code is a fleet code.

\subsection{Window counts and the game for codes}\label{sec:passing-to-codes}

For a fleet code $F$ and $\eps\in E_d$, let
\[
l_F(\eps)=F\cap\bigl(\{0,1,2\}^d+\eps\bigr)\quad\text{and}\quad L(F)=\bigl[\,|l_F(\eps)|:\eps\in E_d\bigr].
\]
The condition $u\in l_F(\eps)$ means $[0,1]^d+u\subseteq[0,3]^d+\eps$, so $L(F)$ is the list of results of the fleet $\ka F_F$ at the origin. Coordinatewise, $u\in l_F(\eps)$ if and only if $u_i\in\{\eps_i,\eps_i+1,\eps_i+2\}$ for all $i$, that is,
\begin{equation}\label{eq:allowed-shifts}
\begin{array}{c|cccc}
u_i&0&1&2&3\\\hline
\text{admissible }\eps_i&\{0\}&\{0,1\}&\{0,1\}&\{1\}
\end{array}
\end{equation}
For fleet codes $F,G$ in the same dimension let
\[
w(F,G)=\bigl|\{(\eps,\eta)\in E_d^2:|l_F(\eps)|>|l_G(\eta)|\}\bigr|\quad\text{and}\quad \Delta(F,G)=w(F,G)-w(G,F).
\]
Thus
\begin{equation}\label{eq:delta-sign}
\Delta(F,G)=\sum_{x\in L(F)}\ \sum_{y\in L(G)}\sgn(x-y),
\end{equation}
where the sums run over the lists with multiplicities. We say that $F$ \emph{beats} $G$ if $\Delta(F,G)>0$. The comparison depends only on the lists $L(F)$ and $L(G)$.

\begin{proposition}\label{prop:count}
For every fleet code $F$ in dimension $d$,
\begin{equation}\label{eq:sum}
\sum_{\eps\in E_d}|l_F(\eps)|=3^d .
\end{equation}
\end{proposition}

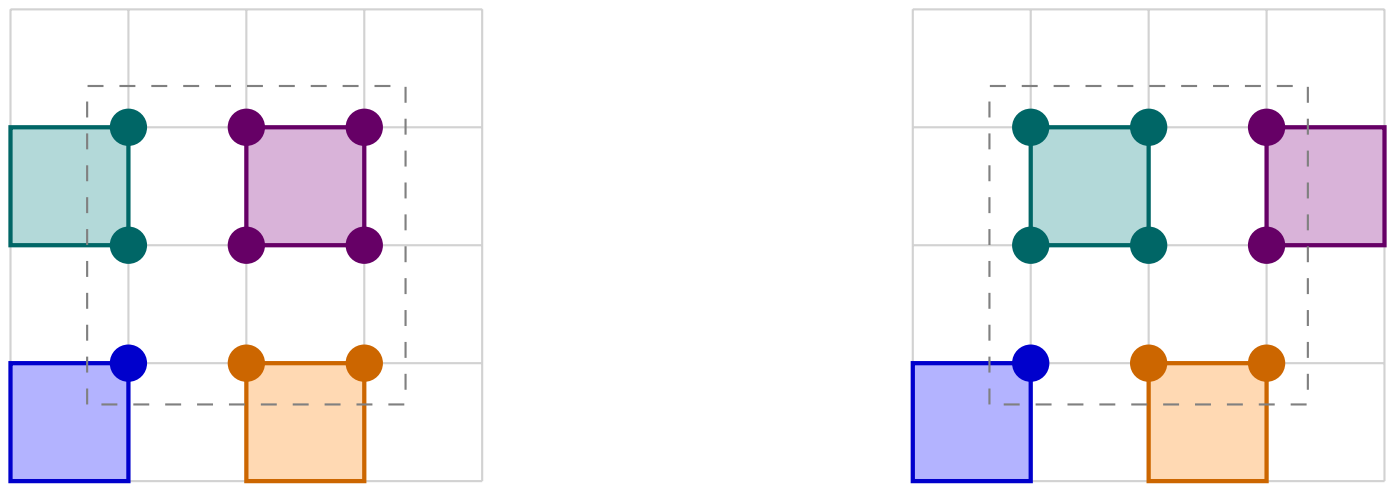
\begin{proof}
By~\eqref{eq:allowed-shifts}, for each fixed letter $u_i$, there are exactly $|I_{u_i}|$ values of $\eps_i\in\{0,1\}$ for which $u_i\in\{\eps_i,\eps_i+1,\eps_i+2\}$. Hence a word $u$ lies in exactly $|I_{u_1}\times\cdots\times I_{u_d}|$ of the sets $l_F(\eps)$. These boxes form a partition of $\{1,2,3\}^d$ by Proposition~\ref{prop:local-code}, so the sum is $3^d$.
\end{proof}

\begin{figure}[ht]
\centering
\localpanel{0/0/blue,2/0/orange,0/2/teal,2/2/violet}{$\ka F^1$: boxes of sizes $1,2,2,4$}
\qquad\qquad
\localpanel{0/0/blue,2/0/orange,1/2/teal,3/2/violet}{$\ka F^2$: boxes of sizes $1,2,4,2$}
\caption{The local codes of Figure~\ref{fig:windows} and the proof of Proposition~\ref{prop:count}. Each of the $9$ lattice points of $\{1,2,3\}^2$ (dashed frame) is a vertex of exactly one black cube. The points of one colour form a box of the minimal partition, and their number equals the number of windows $[0,3]^2+\eps$ containing the corresponding cube. Hence the scores sum to $3^2$.}
\label{fig:local-partition}
\end{figure}

In particular every fleet code has the same mean score,
\begin{equation}\label{eq:mean-score}
\frac1{2^d}\sum_{\eps\in E_d}|l_F(\eps)|=\Bigl(\frac32\Bigr)^d .
\end{equation}
Two further proofs of~\eqref{eq:sum} are given in Sections~\ref{sec:incidence} and~\ref{sec:counting}.

\begin{remark}\label{rem:positive-win-counts}
For any fleet codes $F,G$ in dimension $d\ge1$ we have $w(F,G)>0$. Indeed, the entries of the lists are integers with the non-integral mean $(3/2)^d$, so $\max L(F)>(3/2)^d>\min L(G)$. Thus no code loses every pair of window positions.
\end{remark}

By~\eqref{eq:allowed-shifts}, for a simple code $F=A_1\times\cdots\times A_d$, we have
\[
|l_F(\eps)|=\prod_i c_i(\eps_i),
\]
where $(c_i(0),c_i(1))$ equals $(2,1)$ if $A_i=\{0,2\}$ and $(1,2)$ if $A_i=\{1,3\}$. Hence all simple codes have the list $L(C_d)$, which contains the value $2^j$ exactly $\binom dj$ times (see also Section~\ref{sec:incidence}). In dimension $d=2$, there are exactly two possible lists of results:
\[
L(C_2)=[1,2,2,4]\quad\text{and}\quad L(S)=[1,2,3,3],
\]
where $S=\{00,02,21,23\}$ is the second representative in Proposition~\ref{prop:planar}.

\begin{proposition}\label{prop:extreme-scores}
Let $F$ be a fleet code in dimension $d$.
\begin{enumerate}[label=\textup{(\alph*)}]
\item The list $L(F)$ contains exactly one entry equal to $1$.
\item The list $L(F)$ contains $2^d$ if and only if $F$ is simple, or equivalently if and only if $L(F)=L(C_d)$.
\end{enumerate}
\end{proposition}

\begin{proof}
(a) Let $\mathbf1=(1,\ldots,1)$ and consider the periodic fleet $\ka F_F$ of Proposition~\ref{prop:local-code}. There is a unique cube $Q=[0,1]^d+u$ of this fleet having $2\mathbf1=(2,\ldots,2)$ as a vertex. Then $u\in F\cap\{1,2\}^d$, so $u\in l_F(\eps)$ for every $\eps\in E_d$. Take $\eps=u-\mathbf1$. Every other word $t\in F$ is dichotomous to $u$: in some coordinate, either $u_i=1$ and $t_i=3$, or $u_i=2$ and $t_i=0$. In either case $t\notin l_F(\eps)$, and hence $|l_F(\eps)|=1$. (In other words, the cube $Q$ is the central unit cube of the window $[0,3]^d+\eps$.)

To show that $1$ occurs only once in $L(F)$, suppose that $|l_F(\eps)|=1$, so $Q$ is the only cube of the fleet contained in $[0,3]^d+\eps$. The cubes of the fleet having $\eps+\mathbf1$ and $\eps+2\mathbf1$ as vertices are both contained in this window. They must therefore both be $Q$, which forces $u=\eps+\mathbf1$. Thus $\eps=u-\mathbf1$ is the unique shift with score $1$.

(b) All simple codes have the list $L(C_d)$, whose largest entry is $2^d$. Thus simplicity implies $L(F)=L(C_d)$, which in turn implies $2^d\in L(F)$. Conversely, suppose that $2^d\in L(F)$, and choose $\eps\in E_d$ with $|l_F(\eps)|=2^d$. Then every word of $F$ belongs to $\{0,1,2\}^d+\eps$. Let $V=F-\eps\subseteq\{0,1,2\}^d$. Since distinct words of $F$ are dichotomous, any two distinct words of $V$ have letters $0$ and $2$ in some coordinate. As $|V|=2^d$, the set $V$ is a cube tiling code. Thus $V=\{0,2\}^d$ (see Subsection~\ref{sec:tiling-codes}). Consequently
\[
F=\prod_{i=1}^d\{\eps_i,\eps_i+2\},
\]
that is $F$ is simple.
\end{proof}

\section{Incidence matrices and majorization}\label{sec:majorization}

\subsection{Incidence matrices and cyclic row shifts}\label{sec:incidence}

Order $C_d$ and $E_d$ lexicographically. The \emph{incidence matrix} $M=M_{C_d}$ of the fleet code $C_d$ has rows indexed by $c\in C_d$, columns indexed by $\eps\in E_d$, and $M(c,\eps)=1$ if $c\in l_{C_d}(\eps)$, $M(c,\eps)=0$ otherwise. For $u\in C_d\cup E_d$, let $\supp(u)=\{i:u_i\ne0\}$. By~\eqref{eq:allowed-shifts},
\[
M(c,\eps)=1\iff\supp(\eps)\subseteq\supp(c),
\]
so $M$ is the inclusion matrix of the subsets of $[d]$. The column indexed by $\eps$ has sum $2^{d-|\supp(\eps)|}$; hence $L(C_d)$ contains exactly $\binom dj$ entries equal to $2^j$, $j=0,\dots,d$. For $d=3$, with the column sums in the last row,
\[
\begin{array}{c|rrrrrrrr}
 &000&001&010&011&100&101&110&111\\\hline
000&1&0&0&0&0&0&0&0\\
002&1&1&0&0&0&0&0&0\\
020&1&0&1&0&0&0&0&0\\
022&1&1&1&1&0&0&0&0\\
200&1&0&0&0&1&0&0&0\\
202&1&1&0&0&1&1&0&0\\
220&1&0&1&0&1&0&1&0\\
222&1&1&1&1&1&1&1&1\\\hline
&8&4&4&2&4&2&2&1
\end{array}
\]

Let $z(0)=z(3)=0$, $z(1)=z(2)=2$, and $p(v)=(z(v_1),\ldots,z(v_d))$. For every fleet code $F$ the map $p:F\to C_d$ is a bijection. Indeed, if $u,v\in F$ are distinct, then $u_i=v_i'$ for some $i$, and $z(u_i)\neq z(v_i)$ because $z(0)\ne z(2)$ and $z(1)\ne z(3)$; so $p$ is injective, and $|F|=|C_d|$. The incidence matrix $M_F$ of $F$ is defined in the same way as $M$, with row $v\in F$ placed in the position of row $p(v)$ of $M$. Its column sums are the numbers $|l_F(\eps)|$.

\begin{proposition}\label{prop:shift}
The row of $M_F$ indexed by $v$ is obtained from the row of $M$ indexed by $p(v)$ by a cyclic shift to the right by $k(v)=\sum_{i:\,v_i=3}2^{d-i}$ positions. Equivalently, the binary expansion of $k(v)$ has its ones exactly at the coordinates where $v$ has the letter $3$.
\end{proposition}

\begin{proof}
By~\eqref{eq:allowed-shifts}, both letters $1$ and $2$ allow both possible shifts, $0$ and $1$. Thus replacing $1$ by $2$ in a word does not change the corresponding row of the incidence matrix. Replacing $3$ by $0$ changes the admissible value of $\eps_i$ from $1$ to $0$. Let $\delta_i=1$ if $v_i=3$ and $\delta_i=0$ otherwise, and let $D=\{\eps:M(p(v),\eps)=1\}$. If $\delta_i=1$, then $p(v)_i=0$ and every $\eps\in D$ has $\eps_i=0$. Hence the ones of row $v$ of $M_F$ are in the columns $\eps+\delta$, $\eps\in D$. The columns are arranged in lexicographic order. Number them from $0$ to $2^d-1$. The column labelled by $\eps$ then has index $b(\eps)=\sum_i\eps_i2^{d-i}$. Since $\eps$ and $\delta$ have disjoint supports, $b(\eps+\delta)=b(\eps)+k(v)<2^d$. Thus every one moves exactly $k(v)$ positions to the right without wrapping around, and the entries wrapped around by the cyclic shift are zeros.
\end{proof}

For example, $F=\{000,012,032,120,201,203,222,320\}$ is a fleet code, and with the prescribed order of rows
\[
\begin{array}{c|rrrrrrrr}
 &000&001&010&011&100&101&110&111\\\hline
000&1&0&0&0&0&0&0&0\\
032&0&0&1&1&0&0&0&0\\
320&0&0&0&0&1&0&1&0\\
012&1&1&1&1&0&0&0&0\\
203&0&1&0&0&0&1&0&0\\
201&1&1&0&0&1&1&0&0\\
120&1&0&1&0&1&0&1&0\\
222&1&1&1&1&1&1&1&1\\\hline
&5&4&4&3&4&3&3&1
\end{array}
\]
For instance, $p(032)=002$ and $k(032)=2$; the shifts for $320$ and $203$ are $4$ and $1$, respectively.

\begin{proof}[Second proof of~\eqref{eq:sum}]
By Proposition~\ref{prop:shift}, $M_F$ and $M$ have the same number of ones. Counting the ones of $M$ by columns gives $\sum_{r=0}^d\binom dr2^{d-r}=3^d$.
\end{proof}

\paragraph{From $C_d$ to $F$.}
The passage from $M$ to $M_F$ can be organized one coordinate at a time. For $0\le k\le d$, let
\[
F_k=\{(u_1,\ldots,u_k,z(u_{k+1}),\ldots,z(u_d)):u\in F\}.
\]
Each $F_k$ is a fleet code, since the coordinate maps preserve complementary pairs and hence preserve dichotomy and the number of words. We have $F_0=C_d$ and $F_d=F$. Thus the incidence matrices $M_{F_0},\ldots,M_{F_d}$ start with $M$ and end with $M_F$, restoring the first $k$ coordinates of the words of $F$ at stage $k$.

At step $k$, pair the columns whose indices differ only in coordinate $k$, placing the index with $k$th entry $0$ first. Before this step, every row has one of the forms $(0,0)$, $(1,0)$ or $(1,1)$ on such a pair, since the $k$th letters are $0$ or $2$. Restoring a letter $3$ changes the corresponding pair $(1,0)$ to $(0,1)$, while restoring a letter $1$ leaves the row unchanged. If the two column sums before the step are $A\ge B$, and $t$ of the pairs $(1,0)$ are changed, the sums become
\[
(A,B)\longmapsto(A-t,B+t),\qquad 0\le t\le A-B.
\]
Up to interchanging the two resulting entries, this transfers mass from the larger entry to the smaller one: it is a ``Robin Hood'' transfer~\cite{MarshallOlkinArnold2011}. The pairs of columns are disjoint, so these operations can be performed independently. Proposition~\ref{prop:one-coordinate} below gives a direct partial-sum justification of this averaging effect.

\subsection{The simple code majorizes every fleet code}

Let $x=[x_1,\ldots,x_N]$ and $y=[y_1,\ldots,y_N]$ be real lists in nondecreasing order. We say that $x$ is \emph{majorized} by $y$, and write $x\preccurlyeq y$, if
\[
\sum_{j=1}^sx_j\ge\sum_{j=1}^sy_j\quad(1\le s<N)\quad\text{and}\quad \sum_{j=1}^Nx_j=\sum_{j=1}^Ny_j .
\]
Lists that are not ordered are first arranged in nondecreasing order. This is the usual majorization order~\cite{MarshallOlkinArnold2011}. For fleet codes in the same dimension we write $F\preccurlyeq G$ if $L(F)\preccurlyeq L(G)$. It can easily be shown that the relation $\preccurlyeq$ is a partial order on the equivalence classes of fleet codes in the same dimension, where two codes are equivalent if they have the same list of results.

For a fleet code $F$ and $i\in[d]$, let $F^{(i)}$ be obtained from $F$ by replacing, in the $i$th coordinate of every word, the letter $3$ by $0$ and the letter $1$ by $2$. Both complementary pairs are mapped onto $\{0,2\}$, so $F^{(i)}$ is again a fleet code.

\begin{proposition}\label{prop:one-coordinate}
For every fleet code $F$ and every $i\in[d]$ we have $L(F)\preccurlyeq L(F^{(i)})$.
\end{proposition}

\begin{proof}
For $\eps\in E_d$ with $\eps_i=0$, let $\eps^+$ be obtained by changing $\eps_i$ to $1$. Thus the set $E_d$ is decomposed into pairwise disjoint pairs $\{\eps,\eps^+\}$. On a fixed pair of columns, let $a$, $b$, $c$ count the rows of the forms $(1,0)$, $(1,1)$, $(0,1)$, respectively. The substitutions defining $F^{(i)}$ reverse the operation described in Section~\ref{sec:incidence}: they change $(0,1)$ to $(1,0)$ and leave the other forms unchanged. Thus the column sums change in the following way:
\[
(a+b,b+c)\longmapsto(a+b+c,b).
\]
We perform these operations independently, one pair of columns at a time. After all pairs have been treated, relabelling the rows gives $M_{F^{(i)}}$; the intermediate matrices need not be incidence matrices of fleet codes. For example, when $i=2$, the row labelled $032$ becomes the row labelled $002$ and stays in the same position, since $p(032)=p(002)=002$.

We verify directly that changing one pair cannot increase any sum of the smallest column sums. Let
\[
x_1\le\cdots\le x_N\quad\text{and}\quad y_1\le\cdots\le y_N,\qquad N=2^d,
\]
be the column sums before and after this operation, respectively. All other column sums remain unchanged. Fix $1\le s\le N$. The first $s$ entries correspond to $s$ columns with the smallest sums.

\begin{enumerate}[label=\textup{(\roman*)}]
\item Suppose that neither affected column initially belongs to the first $s$ positions. Then $a+b,b+c\ge x_s$, so $a+b+c$ cannot enter these positions with a value smaller than $x_s$. If $b<x_s$, the new entry $b$ replaces an entry of value $x_s$, and the sum decreases by $x_s-b$. If $b\ge x_s$, the sum is unchanged.

\item Suppose that both affected columns initially belong to the first $s$ positions and both new entries remain there. Their combined contribution is unchanged, since
\[
(a+b)+(b+c)=(a+b+c)+b.
\]

\item Suppose that both affected columns initially belong to the first $s$ positions, but only one new entry remains there. This entry is $b$, while $a+b+c$ leaves. An unchanged entry $z\le a+b+c$ enters in its place, so the change in the sum is
\[
b+z-\bigl((a+b)+(b+c)\bigr)=z-(a+b+c)\le0.
\]

\item Suppose that exactly one affected column initially belongs to the first $s$ positions, say the one with sum $a+b$. Then
\[
b\le a+b\le x_s\le b+c\le a+b+c.
\]
After reordering, $b$ replaces $a+b$ among the first $s$ entries, while $a+b+c$ remains outside them. Thus the sum decreases by $a$. The case where $b+c$ initially belongs to the first $s$ positions is symmetric.
\end{enumerate}

These cases exhaust all possibilities: if both affected columns initially belong to the first $s$ positions, the smaller new entry $b$ can always be retained among them. Hence
\[
\sum_{j=1}^s y_j\le\sum_{j=1}^s x_j.
\]
Applying this inequality successively to all pairs $\{\eps,\eps^+\}$ proves $L(F)\preccurlyeq L(F^{(i)})$. Every operation preserves the total sum, so equality holds for $s=N$.
\end{proof}

\begin{theorem}\label{thm:majorization}
For every fleet code $F$ in dimension $d$,
\[
L(F)\preccurlyeq L(C_d).
\]
That is, if the lists $L(F)=[f_1,\ldots,f_{2^d}]$ and $L(C_d)=[c_1,\ldots,c_{2^d}]$ are given in nondecreasing order, then $\sum_{j\le s}f_j\ge\sum_{j\le s}c_j$ for every $s$, with equality for $s=2^d$.
\end{theorem}

\begin{proof}
Take the codes $F_0,\ldots,F_d$ constructed in Section~\ref{sec:incidence}. For each $1\le k\le d$, applying the replacement in the $k$th coordinate to $F_k$ gives $F_{k-1}$, so Proposition~\ref{prop:one-coordinate} gives
\[
L(C_d)=L(F_0)\succcurlyeq L(F_1)\succcurlyeq\cdots
\succcurlyeq L(F_d)=L(F).
\]
The theorem follows by transitivity.
\end{proof}


\paragraph{Comparable and incomparable codes.}
Majorization is not a total preorder on fleet codes. Let $P=\{000,002,020,022\}$ and
\[
F=P\cup\{201,203,221,223\},\; G=P\cup\{201,203,220,222\}\;\text{and}\; H=P\cup\{200,202,221,223\}.
\]
These are non-simple fleet codes by Proposition~\ref{prop:flip}: $G$ and $H$ arise from $C_3$ by flipping the twin pairs $\{200,202\}$ and $\{220,222\}$, respectively, and $F$ by flipping both. Their lists are
\[
L(F)=[1,2,2,3,3,4,6,6],\qquad L(G)=[1,2,2,3,3,4,5,7]\quad\text{and}\quad L(H)=[1,2,3,3,3,3,5,7].
\]
The partial sums of $L(F)$ and $L(G)$ agree except for $s=7$, where they are $21$ and $20$, respectively; hence $F\preccurlyeq G$. On the other hand, the sums of the first three entries of $L(F)$ and $L(H)$ are $5<6$, and the sums of the first seven are $21>20$, so $F$ and $H$ are incomparable.

\subsection{Consequences of majorization}\label{sec:karamata}

Karamata's inequality states that if $x\preccurlyeq y$ and $\varphi$ is convex on an interval containing the entries of both lists, then $\sum_j\varphi(x_j)\le\sum_j\varphi(y_j)$. For concave $\varphi$ the inequality is reversed. If $\varphi$ is strictly convex or strictly concave, equality holds only if the ordered lists coincide~\cite{KadelburgEtAl2005,MarshallOlkinArnold2011}. By Proposition~\ref{prop:extreme-scores}, equality in the comparison with $L(C_d)$ therefore characterizes simple codes. For a fleet code $F$, we define the random variable $X_F(\eps)=|l_F(\eps)|$, where $\eps$ is chosen uniformly at random from $E_d$.

\begin{corollary}\label{cor:convex-order}
For every fleet code $F$ in dimension $d$ and every convex function $\varphi$ on $(0,\infty)$,
\[
\E[\varphi(X_F)]\le\E[\varphi(X_{C_d})].
\]
For strictly convex $\varphi$, equality holds if and only if $F$ is simple.
\end{corollary}

Applying this inequality to particular convex functions gives the following results.

\begin{corollary}\label{cor:log-mean}
For every fleet code $F$ in dimension $d$,
\[
\mu_F:=\frac1{2^d}\sum_{\eps\in E_d}\ln|l_F(\eps)|\ \ge\ \frac d2\ln2=\mu_{C_d},
\]
with equality if and only if $F$ is simple.
\end{corollary}

\begin{proof}
Apply Corollary~\ref{cor:convex-order} to $\varphi=-\ln$ and use $\sum_{\eps}\ln|l_{C_d}(\eps)|=\sum_k\binom dk\,k\ln2=d2^{d-1}\ln2$.
\end{proof}

\paragraph{Power sums and variance.}
Since $\sum_\eps|l_{C_d}(\eps)|^p=\sum_k\binom dk2^{pk}=(1+2^p)^d$, we get
\[
\sum_{\eps\in E_d}|l_F(\eps)|^p\le(1+2^p)^d\quad(p>1)\quad\text{and}\quad \sum_{\eps\in E_d}|l_F(\eps)|^p\ge(1+2^p)^d\quad(0<p<1),
\]
with equality in both ranges if and only if $F$ is simple. For $p=2$ this gives
\[
\Var(X_F)\le\Bigl(\frac52\Bigr)^d-\Bigl(\frac94\Bigr)^d=\Var(X_{C_d}).
\]
(See also~\cite{DutourItohPoyarkov2007} for the corresponding second-moment bound averaged over all shifts modulo $4$.)
Thus the simple code maximizes the fluctuations of the score in a strong sense. This does not by itself determine the winner of a game; see Sections~\ref{sec:two-coordinates} and~\ref{sec:large-dimensions}. Majorization of lists produced by geometric configurations has also been used, for instance, to compare energies of point configurations on spheres~\cite{Musin2021}.

\section{Codes that beat the simple code}\label{sec:two-coordinates}

For $d\le4$, exhaustive computation (Section~\ref{sec:low-dimensions}) shows that $C_d$ loses to every non-simple fleet code and ties with the simple ones. This fails for $d=5$ (Section~\ref{sec:simple-high}). Nevertheless, in every dimension the codes that differ from $C_d$ in at most two coordinate directions behave as in small dimensions.

\begin{theorem}\label{thm:two-coordinate-changes}
Let $d\ge2$, $C=C_d$, and let $F$ be a fleet code such that, outside two fixed coordinates, every word of $F$ has all its letters in $\{0,2\}$. If $F$ is not simple, then
\[
\Delta(F,C)\ge d-1>0 .
\]
\end{theorem}

The bound is sharp in every dimension (Example~\ref{ex:sharpness}). To prove this inequality, we need a result from extremal combinatorics.

\subsection{A weighted inequality for down-sets}\label{sec:downsets}

A family $\cH\subseteq2^{[r]}$ is a \emph{down-set} if $A\in\cH$ and $B\subseteq A$ imply $B\in\cH$. Write $\binom{[r]}k$ for the family of all $k$-element subsets of $[r]$. For $1\le k\le r$ and $\mathcal G\subseteq\binom{[r]}k$, its \emph{lower shadow} is
\[
\partial\mathcal G=\Bigl\{B\in\binom{[r]}{k-1}:B\subseteq A\text{ for some }A\in\mathcal G\Bigr\}.
\]
The classical local LYM inequality~\cite{Bollobas1986,Engel1997} states that
\[
\frac{|\partial\mathcal G|}{\binom r{k-1}}\ge\frac{|\mathcal G|}{\binom rk}.
\]
If $h_k$ is the number of $k$-element members of a down-set $\cH$, then $\partial(\cH\cap\binom{[r]}k)\subseteq\cH\cap\binom{[r]}{k-1}$, so
\begin{equation}\label{eq:rank-proportions}
\frac{h_0}{\binom r0}\ge\frac{h_1}{\binom r1}\ge\cdots\ge\frac{h_r}{\binom rr}.
\end{equation}

\begin{proposition}\label{prop:downsets}
Let $n=d-2\ge0$, $b_j=\binom dj$ for $0\le j\le d$, and let $\cU\subseteq\cD\subseteq2^{[n]}$ be down-sets with $\emptyset\in\cD$ and $[n]\notin\cU$. Then
\begin{equation}\label{eq:downsets}
\sum_{T\in\cD}(b_{|T|+1}-b_{|T|})+\sum_{T\in\cU}(b_{|T|+2}-b_{|T|+1})\ge d-1.
\end{equation}
\end{proposition}

\begin{proof}
Let $r=n+1=d-1$ and
\[
\cH=\cD\cup\{T\cup\{r\}:T\in\cU\}\subseteq2^{[r]}.
\]
Since $\cU\subseteq\cD$, the family $\cH$ is a down-set. It contains $\emptyset$ but not $[r]$. Let $h_k$ be the number of its $k$-element members and let $c_k=b_{k+1}-b_k$ for $0\le k\le r$. The left-hand side of~\eqref{eq:downsets} equals $\sum_{k=0}^rh_kc_k$.

For $k\le r/2$, inequality~\eqref{eq:rank-proportions} and $\binom rk=\binom r{r-k}$ give $h_k\ge h_{r-k}$. By the symmetry of binomial coefficients, we have $c_{r-k}=-c_k$, and $c_k\ge0$ for $k<r/2$. Since $h_0=1$, $h_r=0$ and $c_0=b_1-b_0=d-1$, pairing the terms of ranks $k$ and $r-k$, we obtain
\[
\sum_{k=0}^rh_kc_k
=\sum_{k=0}^{\lfloor(r-1)/2\rfloor}c_k(h_k-h_{r-k})
\ge c_0(h_0-h_r)=d-1.
\]
If $r$ is even, the middle term vanishes because $c_{r/2}=0$. The same argument applies when $n=0$.
\end{proof}

\subsection{Proof of Theorem~\ref{thm:two-coordinate-changes}}\label{sec:proof-two}

By permuting coordinates, we may assume that the letters $1$ and $3$ occur only in the last two coordinates. Let $n=d-2$ and $b_j=\binom dj$ for $0\le j\le d$. As before, $\{0,2\}^0=\{\two_0\}$.

For $s\in\{0,2\}^n$, let
\[
F_s=\{ab\in\cA^2:sab\in F\}.
\]
Here and below, juxtaposition denotes concatenation of words. The bijection $p:F\to C_d$ of Section~\ref{sec:incidence} fixes the first $n$ coordinates, so $|F_s|=4$. Distinct words of $F_s$ are dichotomous in the last two coordinates. Thus each $F_s$ is a planar fleet code. Conversely, any choice of planar codes $F_s$ gives a fleet code $\{sab:s\in\{0,2\}^n,\ ab\in F_s\}$, and this code is simple if and only if all $F_s$ are equal to one simple planar code.

For $\alpha\beta\in E_2$, let
\[
l_{F_s}(\alpha\beta)=\{ab\in F_s:a\in\{\alpha,\alpha+1,\alpha+2\},\ b\in\{\beta,\beta+1,\beta+2\}\},
\]
and let $B_s$ be the $2\times2$ count matrix with entries $B_s(\alpha,\beta)=|l_{F_s}(\alpha\beta)|$. By Proposition~\ref{prop:planar} and a direct computation, up to interchanging rows, interchanging columns and transposition, $B_s$ has one of the following forms
\begin{equation}\label{eq:two-axes-matrices}
\begin{pmatrix}4&2\\2&1\end{pmatrix}\ \ (F_s\text{ simple})\quad\text{or}\quad
\begin{pmatrix}3&3\\2&1\end{pmatrix}\ \ (F_s\text{ non-simple}).
\end{equation}
Each $B_s$ has exactly one entry $1$. It has an entry $4$ if and only if $F_s$ is simple, and then the position of the $4$ determines $F_s$. Moreover, opposite entries of $B_s$ sum to at least $4$. If the entry $1$ is in the lower right corner, then $B_s$ has one of the following forms
\begin{equation}\label{eq:two-axes-common-minimum}
\begin{pmatrix}4&2\\2&1\end{pmatrix},\qquad
\begin{pmatrix}3&3\\2&1\end{pmatrix}\quad\text{or}\quad
\begin{pmatrix}3&2\\3&1\end{pmatrix}.
\end{equation}

Fix $\delta\in E_n$ and let $T=\{i\in[n]:\delta_i=0\}$. Thus $\delta_i=0$ on $T$ and $\delta_i=1$ outside $T$. If $i\in T$, the window allows both possible letters $s_i=0$ and $s_i=2$ of a prefix $s\in\{0,2\}^n$. If $i\notin T$, then $\delta_i=1$, so the window allows the letters $1,2,3$ in coordinate $i$. Since $s_i\in\{0,2\}$, we must have $s_i=2$. Consequently, the numbers of words of $F$ contained in the four windows indexed by $\delta\alpha\beta$, with $\alpha,\beta\in\{0,1\}$, are the entries of the matrix
\[
B(T)=\sum_{\substack{s\in\{0,2\}^n\\s_i=2\ (i\notin T)}}B_s.
\]
More precisely,
\[
B(T)(\alpha,\beta)=|l_F(\delta\alpha\beta)|,\qquad \alpha\beta\in E_2.
\]

Let $k=|T|$ and $m=2^k$. In the sum defining $B(T)$, each of the $k$ coordinates $s_i$ with $i\in T$ can be chosen independently as $0$ or $2$, while all remaining coordinates are fixed at $2$. Hence the sum has exactly $m$ terms. Let $x_1(T),\ldots,x_4(T)$ be its entries in lexicographic order of the row and column indices. For a fixed $T$, we abbreviate them to $x_1,\ldots,x_4$. Every entry of every $B_s$ belongs to $\{1,2,3,4\}$, so all $x_i$ belong to $\{m,m+1,\ldots,4m\}$. The corresponding four window counts for $C$ are $4m,2m,2m,m$.

Let
\[
\begin{aligned}
\cD&=\{T\subseteq[n]:\min_{1\le i\le4}x_i(T)=2^{|T|}\},\\
\cU&=\{T\subseteq[n]:\max_{1\le i\le4}x_i(T)=4\cdot2^{|T|}\},
\end{aligned}
\]
and let $D_T$ and $U_T$ be the indicator functions of $\cD$ and $\cU$, respectively. The minimum $m$ is attained if and only if all $m$ summands have their entry $1$ in a common position. The maximum $4m$ is attained if and only if all summands are equal to the matrix of one simple planar code. Such a common maximum also gives a common minimum, so $\cU\subseteq\cD$.

If $T'\subseteq T$, the summands of $B(T')$ are among those of $B(T)$. A common position of the entry $1$ or of the entry $4$ is therefore retained when $T$ is replaced by $T'$. Hence both families are down-sets. For $T=\emptyset$, there is only one summand, so $\emptyset\in\cD$. If $[n]\in\cU$, all the sections $F_s$ would be equal to one simple planar code, making $F$ simple. Thus $[n]\notin\cU$. These are precisely the hypotheses of Proposition~\ref{prop:downsets}.

Let $q(T)=\sum_{i=1}^4\sgn(x_i-2m)$. Since opposite entries of $B(T)$ sum to at least $4m$, it follows that $q(T)\ge0$.

If $T\in\cD\setminus\cU$, place the common entry $1$ of the summands in the lower right corner. By~\eqref{eq:two-axes-common-minimum}, the summands have one of the three displayed forms. Since $T\notin\cU$, at least one is of one of the last two types. The upper left entry of $B(T)$ is therefore at least $3m>2m$, the two adjacent entries are at least $2m$ and one of them exceeds $2m$, while the lower right entry is $m<2m$. Hence $q(T)\ge1$. Together with $\cU\subseteq\cD$ and $q(T)\ge0$, this gives
\begin{equation}\label{eq:two-axes-sign-bound}
q(T)\ge D_T-U_T.
\end{equation}

We now express the comparison score in terms of these families. The value $2^j$ occurs $b_j$ times in $L(C)$. Thus a single window count $y>0$ has net score
\[
\gamma_d(y)=\sum_{j=0}^db_j\sgn(y-2^j)
\]
against the entire list $L(C)$: each win contributes $1$, each loss contributes $-1$, and a tie contributes $0$.

For a fixed $T\subseteq[n]$, let $w_T(F)$ be the number of wins of $x_1,x_2,x_3,x_4$ against the entire list $L(C)$, and let $w_T(C)$ be the number of wins of that list against these four counts. In both expressions, $T$ restricts only the window choices for $F$. Their net score is
\[
\Delta_T:=w_T(F)-w_T(C)=\sum_{i=1}^4\gamma_d(x_i).
\]
Every window choice for $F$ belongs to exactly one such group of four, so $\Delta(F,C)=\sum_{T\subseteq[n]}\Delta_T$.

The corresponding four counts for $C$ are $m,2m,2m,4m$, where $m=2^k$ and $k=|T|$. To simplify the contribution of each group, let
\[
a_T=\gamma_d(m)+2\gamma_d(2m)+\gamma_d(4m).
\]
Although $a_T$ need not vanish for an individual $T$, these groups together form the whole list $L(C)$, retaining multiplicities. Consequently, $\sum_{T\subseteq[n]}a_T=\Delta(C,C)=0$: interchanging the two window choices pairs every win with a loss. Hence
\begin{equation}\label{eq:two-axes-centered-score}
\Delta(F,C)=\sum_{T\subseteq[n]}(\Delta_T-a_T).
\end{equation}
Subtracting $a_T$ therefore leaves the total score unchanged. For each $T$, we now compare the contributions of the two groups of four counts against the same entire list $L(C)$.

Expanding $\gamma_d$ gives
\begin{equation}\label{eq:two-axes-expanded-score}
\begin{aligned}
\Delta_T-a_T=\sum_{j=0}^db_j\Bigl[
&\sum_{i=1}^4\sgn(x_i-2^j)-\sgn(m-2^j)\\
&-2\sgn(2m-2^j)-\sgn(4m-2^j)\Bigr].
\end{aligned}
\end{equation}
All counts in both groups belong to $\{m,m+1,\ldots,4m\}$. At a threshold below $m$, their sums of signs are both $4$; at a threshold above $4m$, both are $-4$. These contributions cancel. Thus only the thresholds $m=2^k$, $2m=2^{k+1}$ and $4m=2^{k+2}$ contribute to the sum.

For $p\in\{m,2m,4m\}$, let
\[
\begin{aligned}
s_F(p)&=\sum_{i=1}^4\sgn(x_i-p),\\
s_C(p)&=\sgn(m-p)+2\sgn(2m-p)+\sgn(4m-p).
\end{aligned}
\]
At most one entry of $B(T)$ can equal $m$: equality requires a $1$ in the same position in every summand, and each $B_s$ has a unique entry $1$. Similarly, at most one entry can equal $4m$, since each $B_s$ has at most one entry $4$. Thus the sums of signs at the three thresholds are
\[
\begin{array}{c|c|c|c}
\text{threshold}&s_F(p)&s_C(p)&s_F(p)-s_C(p)\\\hline
p=m&4-D_T&3&1-D_T\\
p=2m&q(T)&0&q(T)\\
p=4m&U_T-4&-3&U_T-1
\end{array}
\]
For example, consider the threshold $p=m$. If $T\in\cD$, exactly one of the four entries $x_i$ equals $m$, while the other three are strictly larger than $m$. One sign is therefore zero and the other three are $1$, giving $s_F(m)=3$. If $T\notin\cD$, all four entries are strictly larger than $m$, so $s_F(m)=4$. Since $D_T$ is the indicator of membership in $\cD$, both cases are expressed by $s_F(m)=4-D_T$. For the simple code, the four counts are $m,2m,2m,4m$, so
\[
s_C(m)=0+1+1+1=3.
\]
Consequently, $s_F(m)-s_C(m)=1-D_T$, as recorded in the first row of the table.

Multiplying the three differences by their respective multiplicities gives
\begin{equation}\label{eq:two-axes-local-score}
\Delta_T-a_T=(1-D_T)b_k+q(T)b_{k+1}-(1-U_T)b_{k+2}.
\end{equation}
Using~\eqref{eq:two-axes-sign-bound} and $b_{k+1}>0$, we obtain
\[
\begin{aligned}
\Delta_T-a_T
&\ge(1-D_T)b_k+(D_T-U_T)b_{k+1}-(1-U_T)b_{k+2}\\
&=b_k-b_{k+2}+D_T(b_{k+1}-b_k)+U_T(b_{k+2}-b_{k+1}).
\end{aligned}
\]
We now sum over all $T\subseteq[n]$. The terms independent of $D_T$ and $U_T$ cancel, since
\[
\sum_{T\subseteq[n]}(b_{|T|}-b_{|T|+2})=\sum_{k=0}^n\binom nk(b_k-b_{k+2})=0.
\]
Indeed, the terms for $k$ and $n-k$ have equal weights and opposite signs, because $d=n+2$ and the binomial coefficients are symmetric. The indicators $D_T$ and $U_T$ select precisely the sets in $\cD$ and $\cU$, respectively. Together with~\eqref{eq:two-axes-centered-score}, this gives
\begin{equation}\label{eq:two-axes-score-bound}
\Delta(F,C)\ge\sum_{T\in\cD}(b_{|T|+1}-b_{|T|})+\sum_{T\in\cU}(b_{|T|+2}-b_{|T|+1}).
\end{equation}
The hypotheses of Proposition~\ref{prop:downsets} were verified above, so~\eqref{eq:downsets} bounds the right-hand side by $d-1$ from below. This completes the proof. \qed

\medskip
The following construction shows the sharpness of the bound $d-1$ in Theorem~\ref{thm:two-coordinate-changes} in every dimension $d\ge2$.

\begin{example}\label{ex:sharpness}
Let $d\ge2$ and
\[
F=\bigl(C_d\setminus\{\two_{d-1}0,\two_{d-1}2\}\bigr)\cup\{\two_{d-1}1,\two_{d-1}3\}.
\]
This is a non-simple fleet code whose letters $1,3$ occur only in the last coordinate. For $\delta\in E_{d-1}$, let $r=d-1-|\delta|$. The counts in the windows $\delta0,\delta1$ change from $2^{r+1},2^r$ to $2^{r+1}-1,2^r+1$, respectively. If $r=0$, the two counts are interchanged, so their total contribution to the comparison with $C_d$ is unchanged. If $r\ge1$, only comparisons with counts equal to $2^{r+1}$ or $2^r$ change. Since the value $2^j$ occurs $\binom dj$ times in $L(C_d)$, we obtain
\[
\Delta(F,C_d)=\sum_{r=1}^{d-1}\binom{d-1}{r}\left(\binom dr-\binom d{r+1}\right)=d-1.
\]
Indeed, the full sum over $0\le r\le d-1$ vanishes by pairing $r$ with $d-1-r$, and the omitted term is $1-d$.
\end{example}

\begin{remark}
The role of dimension $2$ in Theorem~\ref{thm:two-coordinate-changes} is visible in~\eqref{eq:two-axes-matrices} and~\eqref{eq:two-axes-common-minimum}: these provide the diagonal-sum bound $q(T)\ge0$ and the stronger bound $q(T)\ge1$ for $T\in\cD\setminus\cU$. Writing higher-dimensional sections as arrays does not supply the same inequalities. Indeed, the analogue of Theorem~\ref{thm:two-coordinate-changes} for changes in three coordinates fails: the code $H_5$ of Section~\ref{sec:simple-high} uses the letters $1$ and $3$ only in the coordinates $1$, $4$ and $5$, yet $C_5$ beats it.
\end{remark}

\section{Exhaustive comparisons in dimensions at most four}\label{sec:low-dimensions}

Recall that $\mathscr F_d$ is the set of fleet codes in dimension $d$ with fixed coordinate and letter labels. Let $N_d=|\mathscr F_d|$. A code $F$ is \emph{unbeaten} if $\Delta(F,G)\ge0$ for every $G\in\mathscr F_d$. Let
\[
\rho_d(F)=\frac{|\{G\in\mathscr F_d\setminus\{F\}:\Delta(F,G)>0\}|}{N_d-1}
\]
be the proportion of opponents beaten by $F$; tied opponents remain in the denominator.

To compute the game comparisons, we enumerated $\mathscr F_d$ for $d\le4$. 
Codes with equal lists have equal scores against every opponent, so for $d=4$ the comparisons reduce to $2969$ distinct lists, with their multiplicities retained. 
\begin{center}
\begin{tabular}{crrrc}
\toprule
$d$ & $N_d$ & distinct lists & unbeaten codes & $\min\Delta(F,C_d)$, $F$ non-simple\\
\midrule
1 & 2 & 1 & 2 & ---\\
2 & 12 & 2 & 8 & 1\\
3 & 744 & 14 & 0 & 2\\
4 & $5\,541\,744$ & 2969 & 0 & 3\\
\bottomrule
\end{tabular}
\captionof{table}{Enumeration and comparison results for $d\le4$.}
\end{center}

\paragraph{Comparison with the simple code.}
For $2\le d\le4$, every non-simple $F\in\mathscr F_d$ satisfies $\Delta(F,C_d)\ge d-1$, not only the codes covered by Theorem~\ref{thm:two-coordinate-changes}. 
By Proposition~\ref{prop:extreme-scores}, the codes $F$ with $L(F)=L(C_d)$ are exactly the $2^d$ simple codes in every dimension. In particular, $C_4$ loses to all $5\,541\,728$ non-simple codes and ties with the other $15$ simple codes.

\paragraph{Unbeaten codes and the largest winning proportion.}
For $d=1$ all comparisons are tied. For $d=2$ the eight non-simple codes, all with the list $[1,2,3,3]$, beat the four simple codes and tie with each other. For $d=3$ and $d=4$ every code is beaten by at least $48$ and $59\,520$ codes, respectively, so no unbeaten code exists. For $d=4$,
\[
\max_{F\in\mathscr F_4}\rho_4(F)=\frac{5\,317\,104}{5\,541\,743}\approx95.95\%,
\]
attained by exactly $3072$ codes, all with the list $[1,2,2,3,4,4,4,5,6,6,6,7,7,8,8,8]$. One of them is
\begin{equation}\label{eq:maximizer-V}
\begin{aligned}
V=\{&3213,3033,2321,1313,1301,3201,0233,1033,\\
    &2233,0121,1113,3013,3001,2121,1101,0321\},
\end{aligned}
\end{equation}
which beats $5\,317\,104$ codes, loses to $59\,520$ and ties with $165\,119$. 

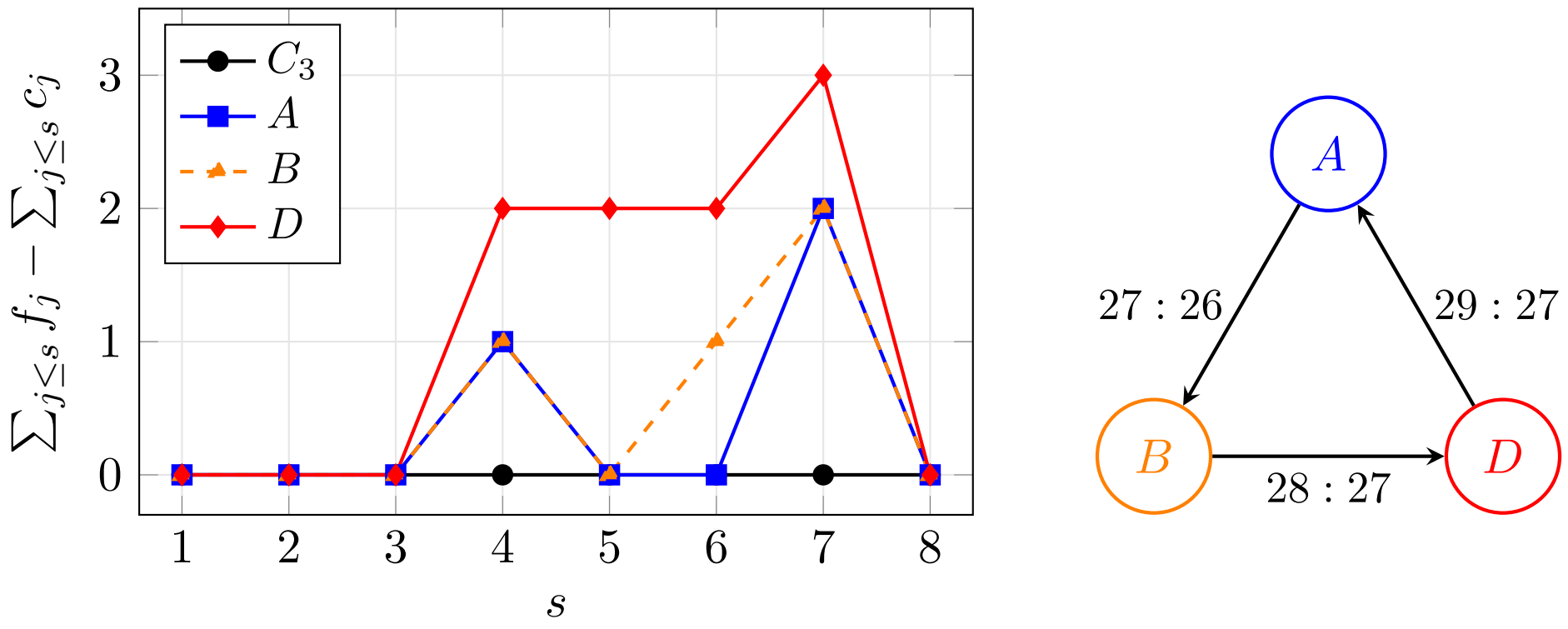
\begin{figure}[ht]
\centering
\begin{tikzpicture}
\begin{axis}[width=8.2cm,height=5.6cm,
  xlabel={$s$}, ylabel={$\sum_{j\le s}f_j-\sum_{j\le s}c_j$},
  xmin=0.6,xmax=8.4,ymin=-0.3,ymax=3.5,xtick={1,...,8},ytick={0,1,2,3},
  grid=major, grid style={gray!20}, legend pos=north west, legend cell align=left,
  legend style={font=\small}]
\addplot[black,mark=*,thick] coordinates {(1,0)(2,0)(3,0)(4,0)(5,0)(6,0)(7,0)(8,0)};
\addplot[blue,mark=square*,thick] coordinates {(1,0)(2,0)(3,0)(4,1)(5,0)(6,0)(7,2)(8,0)};
\addplot[orange,mark=triangle*,thick,dashed] coordinates {(1,0)(2,0)(3,0)(4,1)(5,0)(6,1)(7,2)(8,0)};
\addplot[red,mark=diamond*,thick] coordinates {(1,0)(2,0)(3,0)(4,2)(5,2)(6,2)(7,3)(8,0)};
\legend{$C_3$,$A$,$B$,$D$}
\end{axis}
\end{tikzpicture}
\hspace{0.6cm}
\begin{tikzpicture}[baseline=-2.2cm, >=stealth]
  \node[circle,draw,thick,blue,minimum size=0.9cm] (A) at (90:1.6) {$A$};
  \node[circle,draw,thick,orange,minimum size=0.9cm] (B) at (210:1.6) {$B$};
  \node[circle,draw,thick,red,minimum size=0.9cm] (D) at (330:1.6) {$D$};
  \draw[->,thick] (A) -- node[left,font=\small]{$27:26$} (B);
  \draw[->,thick] (B) -- node[below,font=\small]{$28:27$} (D);
  \draw[->,thick] (D) -- node[right,font=\small]{$29:27$} (A);
\end{tikzpicture}
\caption{Left: excess of the partial sums of $L(A)$, $L(B)$, $L(D)$ over those of $L(C_3)=[1,2,2,2,4,4,4,8]$. The curves are nested, so $D\preccurlyeq B\preccurlyeq A\preccurlyeq C_3$. Right: nevertheless the three codes form a cycle; an arrow $X\to Y$ means that $X$ beats $Y$, labelled $w(X,Y):w(Y,X)$. 
}
\label{fig:cycle}
\end{figure}

\paragraph{Nontransitivity.}
Let us consider the following three fleet codes:
\[
\begin{aligned}
A&=\{013,033,111,131,213,233,311,331\},\\
B&=\{103,111,123,131,311,313,330,332\},\\
D&=\{013,033,111,131,203,223,311,331\},
\end{aligned}
\]
with lists $L(A)=[1,2,2,3,3,4,6,6]$, $L(B)=[1,2,2,3,3,5,5,6]$ and $L(D)=[1,2,2,4,4,4,5,5]$. Then
\[
\Delta(A,B)=27-26=1,\qquad \Delta(B,D)=28-27=1\quad\text{and}\quad \Delta(D,A)=29-27=2,
\]
so $A$ beats $B$, $B$ beats $D$, and $D$ beats $A$. The partial sums of $L(A)$, $L(B)$ and $L(D)$ are, respectively,
\[
[1,3,5,8,11,15,21,27],\qquad [1,3,5,8,11,16,21,27]\quad\text{and}\quad [1,3,5,9,13,17,22,27].
\]
Hence $D\preccurlyeq B\preccurlyeq A\preccurlyeq C_3$, with all lists distinct: the cycle runs along a chain of the majorization order (Figure~\ref{fig:cycle}). The code $D$ attains $\max\rho_3=624/743\approx83.98\%$, with $48$ losses and $71$ ties. There is also a cycle through the $4$-dimensional maximizer $V$ in~\eqref{eq:maximizer-V}: for
\[
\begin{aligned}
U=\{&0022,0222,1000,1020,1102,1210,1230,1302,\\
    &2022,2222,3030,3102,3110,3230,3302,3310\},\\
W=\{&0022,0101,0120,0222,0330,1103,1302,1310,\\
    &2022,2101,2120,2222,2330,3103,3302,3310\}
\end{aligned}
\]
we have $\Delta(U,V)=\Delta(V,W)=\Delta(W,U)=1$.

\paragraph{Dependence on the isomorphic copy.}
Let
\[
F=\{013,111,131,133,213,311,331,333\}\quad\text{and}\quad G=\{011,113,131,133,211,313,331,333\}.
\]
Interchanging $1$ and $3$ in the third coordinate maps $F$ onto $G$, so the codes are isomorphic. Their lists are $[1,2,2,3,3,4,5,7]$ and $[1,2,3,3,3,3,5,7]$, respectively, and against the code $B$ above we have
\[
\Delta(F,B)=-1\quad\text{and}\quad \Delta(G,B)=1 .
\]
Thus comparisons are not invariant under isomorphism.

\paragraph{Comparability.}
Let $p_d$ be the probability that two distinct codes chosen uniformly from $\mathscr F_d$ are comparable in the majorization order, pairs with identical lists included. Computing over all codes, we obtain $p_2=1$,
\[
p_3=\frac{30\,892+222\,464}{276\,396}\approx0.9166\quad\text{and}\quad 
p_4=\frac{17\,704\,925\,896+9\,639\,138\,755\,072}{15\,355\,460\,509\,896}\approx0.6289,
\]
where the two summands count the pairs with identical lists and with distinct comparable lists, respectively. To compute $p_4$, it suffices to compare the $2969$ distinct lists: a comparable pair of distinct lists with multiplicities $m_i,m_j$ contributes $m_i m_j$ pairs of codes, while a list of multiplicity $m_i$ contributes $\binom{m_i}{2}$ pairs with identical lists. The values for $d\le4$ suggest that comparability becomes less frequent as the dimension grows; see Problem~\ref{prob:comparability}.

\subsection{A geometric realization of the Steinhaus--Trybu{\l}a paradox}\label{sec:geometric-paradox}

The Steinhaus--Trybu{\l}a paradox~\cite{SteinhausTrybula1959,Trybula1961} is the existence of independent random variables $X_1,X_2,X_3$ such that
\[
\PP(X_1>X_2)>\frac12,\qquad
\PP(X_2>X_3)>\frac12\quad\text{and}\quad 
\PP(X_3>X_1)>\frac12.
\]
These inequalities are stronger than the cyclic comparisons defined by $w(F,G)>w(G,F)$ when draws are possible. The cycle above in dimension $d=3$ has winning probabilities $27/64$, $28/64$ and $29/64$. The following construction realizes the stronger inequalities by geometric window counts.

\begin{proposition}\label{prop:geometric-paradox}
There exist three $4\Z^4$-periodic fleets $\ka F_1,\ka F_2,\ka F_3$ in $\R^4$ such that, for independent uniform shifts $\eps^1,\eps^2,\eps^3\in E_4$ and scores $X_i(\eps^i)=N_0(\ka F_i,\eps^i)$,
\begin{equation}\label{eq:geometric-paradox}
\begin{aligned}
\PP(X_1>X_2)&=\frac{132}{256}>\frac12,\\
\PP(X_2>X_3)&=\frac{131}{256}>\frac12,\\
\PP(X_3>X_1)&=\frac{129}{256}>\frac12.
\end{aligned}
\end{equation}
In particular, all three probabilities are greater than $1/2$, while $\E[X_i]=81/16$ for each $i$.
\end{proposition}

\begin{proof}
Consider the following fleet codes:
\[
\begin{aligned}
F_1=\{&0000,0002,0020,0022,0200,0232,3220,0212,\\
    &2300,2332,2020,2012,2100,2132,1220,2212\},\\[2pt]
F_2=\{&0000,0002,0020,0022,0203,3201,3223,0221,\\
    &2303,2001,2023,2321,2103,1201,1223,2121\},\\[2pt]
F_3=\{&0000,0002,0023,0021,0233,3231,3213,0211,\\
    &2333,2031,2010,2012,2133,1231,1213,2211\}.
\end{aligned}
\]
For each $i\in\{1,2,3\}$, Proposition~\ref{prop:local-code} defines the periodic fleet
\[
\ka F_i=\{[0,1]^4+u+4z:\ u\in F_i,\ z\in\Z^4\}.
\]
Counting the words of $F_i$ contained in $\{0,1,2\}^4+\eps$ for the $16$ shifts $\eps\in E_4$ gives the following lists, displayed in nondecreasing order:
\[
\begin{aligned}
L(F_1)&=[1,2,2,3,3,5,5,5,5,6,6,6,6,6,9,11],\\
L(F_2)&=[1,2,3,4,4,4,4,4,4,4,5,8,8,8,8,10],\\
L(F_3)&=[1,2,3,3,3,3,3,4,6,7,7,7,7,7,8,10].
\end{aligned}
\]
Every list has sum $81$. Comparing their entries gives
\begin{center}
\begin{tabular}{c@{\qquad}rrr}
\toprule
Comparison & Wins & Losses & Draws\\
\midrule
$F_1$ against $F_2$ & 132 & 115 & 9\\
$F_2$ against $F_3$ & 131 & 106 & 19\\
$F_3$ against $F_1$ & 129 & 109 & 18\\
\bottomrule
\end{tabular}
\end{center}
Dividing the numbers of wins by $256$ proves the inequalities~\eqref{eq:geometric-paradox}.
\end{proof}

In the classical dice interpretation, the random variables represent the outcomes of rolls. Here their values arise from counting black cubes in translated windows, so the cyclic inequalities are determined by the geometry of the fleets.


\section{Large dimensions}\label{sec:large-dimensions}

\subsection{Cartesian powers}\label{sec:powers}

For a fleet code $F$ in dimension $k$ and $n\ge1$ let 
$$F^n=\{u^1\cdots u^n:u^1,\ldots,u^n\in F\}\subseteq\cA^{kn}.$$
Two distinct words of $F^n$ differ in some block, where they are dichotomous because $F$ is a fleet code, so $F^n$ is a fleet code. By~\eqref{eq:allowed-shifts}, a word belongs to a window if and only if each of its blocks belongs to the corresponding window, so
\[
l_{F^n}(\eps^1\cdots\eps^n)=l_F(\eps^1)\times\cdots\times l_F(\eps^n)\quad\text{and}\quad |l_{F^n}(\eps^1\cdots\eps^n)|=\prod_{i=1}^n|l_F(\eps^i)| .
\]
Thus the score of $F^n$ is the product of $n$ independent scores of $F$, each obtained by choosing a shift uniformly from $E_k$. Note that $C_k^n=C_{kn}$. Recall that
\[
\mu_F=\frac{1}{2^k}\sum_{\eps\in E_k}\ln|l_F(\eps)|
\]
is the mean logarithm of the score of $F$.

By taking Cartesian powers, we can obtain pairs of fleets for which the probability that one loses to the other tends to $1$ as the dimension grows. The following theorem makes this precise in terms of the mean logarithms of their scores.

\begin{theorem}\label{thm:lln-products}
Let $F,G$ be fleet codes in dimension $k$ with $\mu_F<\mu_G$. Then
\[
\lim_{n\to\infty}\frac{w(G^n,F^n)}{4^{kn}}=1,\qquad \lim_{n\to\infty}\frac{w(F^n,G^n)}{4^{kn}}=0\quad\text{and}\quad \lim_{n\to\infty}\frac{w(G^n,F^n)}{w(F^n,G^n)}=+\infty .
\]
In particular, the proportion of tied pairs tends to $0$.
\end{theorem}

\begin{proof}
Let $\eps^i,\delta^i\in E_k$, $i\ge1$, be mutually independent shifts, each chosen uniformly at random. Define the random variables $S_i(\eps^i)=|l_F(\eps^i)|$ and $P_i(\delta^i)=|l_G(\delta^i)|$. Then $\prod_{i\le n}S_i$ and $\prod_{i\le n}P_i$ have the score distributions of $F^n$ and $G^n$, respectively, under independent uniform shifts. Define $Z_i(\eps^i,\delta^i)=\ln P_i(\delta^i)-\ln S_i(\eps^i)$. These variables are independent and identically distributed, $|Z_i|\le k\ln2$ because all scores lie in $\{1,...,2^k\}$, and $\E[Z_i]=\mu_G-\mu_F=:\theta>0$. By the weak law of large numbers,
\[
\PP\Bigl(\prod_{i\le n}P_i\le\prod_{i\le n}S_i\Bigr)=\PP\Bigl(\frac1n\sum_{i\le n}Z_i\le0\Bigr)\le\PP\Bigl(\Bigl|\frac1n\sum_{i\le n}Z_i-\theta\Bigr|\ge\theta\Bigr)\longrightarrow0 .
\]
This gives the first two limits. The third follows, since $w(F^n,G^n)>0$ by Remark~\ref{rem:positive-win-counts}.
\end{proof}

For example, for $C_2=\{00,02,20,22\}$ and $S=\{00,02,21,23\}$, with lists $[1,2,2,4]$ and $[1,2,3,3]$, respectively, we have $\mu_S-\mu_{C_2}=\frac14\ln\frac{18}{16}>0$. Hence $S^n$ beats $C_{2n}$ with probability tending to $1$.

Thus sufficiently high Cartesian powers eliminate nontransitive cycles among fleet codes with pairwise distinct mean logarithms of the scores:

\begin{corollary}\label{cor:powers-order}
Let $F_1,\ldots,F_q$ be fleet codes in the same dimension with pairwise distinct values $\mu_{F_1}>\cdots>\mu_{F_q}$. There is $n_0$ such that for all $n\ge n_0$ and all $i<j$, the code $F_i^n$ beats $F_j^n$. In particular, sufficiently high powers of a nontransitive cycle with distinct values of $\mu$ are no longer cyclic.
\end{corollary}

\begin{proof}
Apply Theorem~\ref{thm:lln-products} to each of the finitely many pairs.
\end{proof}

For example, for the cycle $A,B,D$ of Section~\ref{sec:low-dimensions}, the products of the entries of the lists are $5184$, $5400$ and $6400$, respectively, so $\mu_A<\mu_B<\mu_D$. Hence for large $n$ the code $D^n$ beats $B^n$ and $A^n$, and $B^n$ beats $A^n$. Here the strict order of the means follows from the majorization chain $D\preccurlyeq B\preccurlyeq A$ and the strict concavity of the logarithm.

\begin{corollary}\label{cor:powers-simple}
Let $F$ be a fleet code in dimension $d$.
If $F$ is not simple, then $F^n$ beats $C_{dn}$ with probability tending to $1$. In particular, if $C_d$ beats $F$, then the comparison is reversed for all sufficiently high powers.
\end{corollary}

\begin{proof}
This follows from Corollary~\ref{cor:log-mean} and Theorem~\ref{thm:lln-products}.
\end{proof}

Thus the simple code is asymptotically the weakest under Cartesian powers: every fixed non-simple code eventually beats it with probability tending to one, even if the initial comparison favours the simple code. Theorem~\ref{thm:simple-asymptotic}, which we prove in the next subsection, complements this observation by allowing the opposing code to vary with the dimension: the ratio of wins to losses of the simple code has limit superior at most one.

\subsection{The simple code in high dimensions}\label{sec:simple-high}

Theorem~\ref{thm:two-coordinate-changes} does not extend to arbitrary non-simple codes. In dimension $d=5$ let
\begin{equation}\label{eq:A5-B5}
\begin{aligned}
A_5&=\{22000,22002,00022,20022,20200,20220,02222,22222\},\\
B_5&=\{22003,22001,30022,10022,20230,20210,32222,12222\},
\end{aligned}
\end{equation}
and $H_5=(C_5\setminus A_5)\cup B_5$. The set $A_5$ consists of four twin pairs of $C_5$, and $B_5$ is obtained by flipping them, so $H_5$ is a fleet code by Proposition~\ref{prop:flip}. Its list is
\[
L(H_5)=[1,2,3,3,3,3,3,3,3,3,5,5,5,6,6,6,6,7,7,7,8,8,8,8,11,11,14,14,15,15,16,28],
\]
and $w(H_5,C_5)=486<487=w(C_5,H_5)$. Thus $C_5$ beats $H_5$. This advantage disappears for powers: since
\[
\mu_{H_5}-\mu_{C_5}=\frac{14\ln3+5\ln5+6\ln7+2\ln11-55\ln2}{32}\approx0.0555>0,
\]
Theorem~\ref{thm:lln-products} gives $w(H_5^n,C_{5n})/w(C_{5n},H_5^n)\to+\infty$.

The following theorem holds for arbitrary sequences of codes: although $C_d$ can beat an opponent, the ratio of its wins to its losses cannot stay above $1$ by a fixed amount.

\begin{theorem}\label{thm:simple-asymptotic}
For every sequence of fleet codes $F_d\in\mathscr F_d$,
\[
\limsup_{d\to\infty}\frac{w(C_d,F_d)}{w(F_d,C_d)}\le1 .
\]
More precisely, for every fleet code $F$ in dimension $d\ge1$ and every integer $1\le m\le d$,
\begin{equation}\label{eq:uniform-bound}
\frac{w(F,C_d)}{4^d}\ \ge\ \frac12-2^{-m}-\Bigl(m+\frac12\Bigr)\binom{2d}d4^{-d}.
\end{equation}
\end{theorem}

\begin{proof}
Let $N=2^d$, $h=2^{d-m}$, and let the lists $L(F)=[f_1,\ldots,f_N]$ and $L(C_d)=[c_1,\ldots,c_N]$ be arranged in nondecreasing order. For $1\le j\le N-h$, Theorem~\ref{thm:majorization} gives
\[
f_{j+h}\ \ge\ \frac1{j+h}\sum_{i\le j+h}f_i\ \ge\ \frac1{j+h}\sum_{i\le j+h}c_i\ \ge\ \frac h{j+h}\,c_j\ \ge\ \frac hNc_j=2^{-m}c_j,
\]
where the third inequality uses $c_i\ge c_j$ for $j<i\le j+h$. Hence $c_j>2^mc_\ell$ implies $f_{j+h}>c_\ell$. Discarding the $hN$ pairs $(j,\ell)$ with $j>N-h$, we get
\begin{equation}\label{eq:far-wins}
w(F,C_d)\ge\bigl|\{(j,\ell):c_j>2^mc_\ell\}\bigr|-hN .
\end{equation}
The value $2^r$ occurs $\binom dr$ times in $L(C_d)$, so the number of pairs with $c_j=2^tc_\ell$ is 

$$T_t=\sum_r\binom dr\binom d{r+t}=\binom{2d}{d+t}\le\binom{2d}d=T_0$$
 by Vandermonde's identity. Exactly $(N^2-T_0)/2$ pairs satisfy $c_j>c_\ell$, so
\[
\bigl|\{(j,\ell):c_j>2^mc_\ell\}\bigr|=\frac{N^2-T_0}2-\sum_{t=1}^mT_t\ \ge\ \frac{N^2}2-\Bigl(m+\frac12\Bigr)T_0 .
\]
Together with~\eqref{eq:far-wins} and $h/N=2^{-m}$, this proves~\eqref{eq:uniform-bound}.

\smallskip
Since $\binom{2d}d4^{-d}\to0$, by \eqref{eq:uniform-bound} we have 
$$\liminf_dq_d\ge\frac12-2^{-m}$$ 
for every $m$, where $q_d=w(F_d,C_d)/4^d$. Hence $\liminf_dq_d\ge\frac12$. As $w(C_d,F_d)\le4^d-w(F_d,C_d)$ and $q_d>0$,
\[
\limsup_{d\to\infty}\frac{w(C_d,F_d)}{w(F_d,C_d)}\le\limsup_{d\to\infty}\Bigl(\frac1{q_d}-1\Bigr)\le1 . \qedhere
\]
\end{proof}

Theorem~\ref{thm:simple-asymptotic} does not bound the unnormalized advantage $\Delta(C_d,F_d)$. Applying the four flips that define $H_5$ to a suitable set of words in $C_d$ shows that this advantage can be arbitrarily large:

\begin{proposition}\label{prop:unbounded-simple-advantage}
For every $d\ge5$ there exists a fleet code $F_d\in\mathscr F_d$ such that
\begin{equation}\label{eq:appendix-exact-advantage}
\Delta(C_d,F_d)=\binom{2d-5}{d-4}-d+1 .
\end{equation}
Thus $\Delta(C_d,F_d)\to\infty$, while $\Delta(C_d,F_d)/4^d\to0$, in accordance with Theorem~\ref{thm:simple-asymptotic}.
\end{proposition}

\begin{proof}
We begin with the counts for $C_5$ and $H_5$ at each shift $\eta\in E_5$, listed in Table~\ref{tab:C5-H5-counts}.
\begin{center}
\small
\begin{tabular}{crr@{\qquad\qquad}crr}
\toprule
$\eta$ & $|l_{C_5}(\eta)|$ & $|l_{H_5}(\eta)|$ & $\eta$ & $|l_{C_5}(\eta)|$ & $|l_{H_5}(\eta)|$\\
\midrule
$00000$ & 32 & 28 & $10000$ & 16 & 16\\
$00001$ & 16 & 15 & $10001$ & 8 & 11\\
$00010$ & 16 & 15 & $10010$ & 8 & 11\\
$00011$ & 8 & 6 & $10011$ & 4 & 6\\
$00100$ & 16 & 14 & $10100$ & 8 & 8\\
$00101$ & 8 & 7 & $10101$ & 4 & 5\\
$00110$ & 8 & 8 & $10110$ & 4 & 6\\
$00111$ & 4 & 3 & $10111$ & 2 & 3\\
$01000$ & 16 & 14 & $11000$ & 8 & 8\\
$01001$ & 8 & 8 & $11001$ & 4 & 6\\
$01010$ & 8 & 7 & $11010$ & 4 & 5\\
$01011$ & 4 & 3 & $11011$ & 2 & 3\\
$01100$ & 8 & 7 & $11100$ & 4 & 5\\
$01101$ & 4 & 3 & $11101$ & 2 & 3\\
$01110$ & 4 & 3 & $11110$ & 2 & 3\\
$01111$ & 2 & 1 & $11111$ & 1 & 2\\
\bottomrule
\end{tabular}
\captionof{table}{Counts in lexicographic order of $\eta$, read down the left block and then down the right block.}
\label{tab:C5-H5-counts}
\end{center}

Let $n=d-5$ and, using the sets in~\eqref{eq:A5-B5}, define
\[
F_d=\bigl(C_d\setminus(\{\two_n\}\times A_5)\bigr)\cup\bigl(\{\two_n\}\times B_5\bigr).
\]
Thus we perform the same four flips as for $H_5$, with the prefix $\two_n$ attached to every word involved. Each flipped pair remains a twin pair, so $F_d$ is a fleet code by Proposition~\ref{prop:flip}.

We first compare the counts of $C_d$ and $F_d$ for the same window shift $\eps=\delta\eta$, where $\delta\in E_n$ and $\eta\in E_5$. Every changed word has the prefix $\two_n$, whose letters are admissible for either shift, $0$ or $1$. Consequently,
\[
|l_{F_d}(\delta\eta)|-|l_{C_d}(\delta\eta)|
=|l_{H_5}(\eta)|-|l_{C_5}(\eta)|=:e(\eta).
\]
Thus the difference depends only on $\eta$, although the counts themselves also depend on $\delta$. A decrease or increase below always refers to the change from $C_d$ to $F_d$ at this same shift.

Let
\[
Z=\{00110,01001,10000,10100,11000\}.
\]
Table~\ref{tab:C5-H5-counts} and the identity defining $e(\eta)$ show that, for every $\delta\in E_n$,
\[
\begin{aligned}
|l_{C_d}(\delta\eta)|&>|l_{F_d}(\delta\eta)|
&&\text{if }\eta\notin Z\text{ and }\eta_1=0,\\
|l_{C_d}(\delta\eta)|&<|l_{F_d}(\delta\eta)|
&&\text{if }\eta\notin Z\text{ and }\eta_1=1,\\
|l_{C_d}(\delta\eta)|&=|l_{F_d}(\delta\eta)|
&&\text{if }\eta\in Z.
\end{aligned}
\]
Equivalently, $e(\eta)=0$ on $Z$ and $\sgn(e(\eta))=2\eta_1-1$ elsewhere.

If $\delta\in E_n$ has $k$ zero coordinates and $s=|\eta|$, the count in the window $\delta\eta$ changes from $2^{k+5-s}$ to
\[
|l_{F_d}(\delta\eta)|=2^{k+5-s}+e(\eta).
\]
For $\eta\notin Z$, the differences in Table~\ref{tab:C5-H5-counts} show that the new count lies strictly between the original power of two and the adjacent power in the direction of the change, except when $k=0$ and $\eta\in\{01111,11111\}$. These two exceptional counts change from $2$ to $1$ and from $1$ to $2$, respectively. To see that there are no further exceptions, we note that, for fixed $\eta$, the correction $e(\eta)$ is independent of $k$, whereas the distances from $2^{k+5-s}$ to the adjacent powers of two double whenever $k$ increases by one. Thus every strict inequality checked at $k=0$ remains strict for larger $k$, and the two exceptional equalities at $k=0$ become strict inequalities for every $k\ge1$. We account for the two exceptions at the end.

We compute $\Delta(C_d,F_d)$ by starting from $\Delta(C_d,C_d)=0$ and changing only the second list. If a count $2^a$ decreases without reaching $2^{a-1}$, its $\binom da$ draws against $C_d$ become wins for $C_d$, and all other comparisons stay unchanged. An increase that does not reach $2^{a+1}$ turns the same draws into losses. For a fixed $\eta$ of weight $s$, summing this rule over the additional shifts $\delta$ therefore gives the contribution
\[
-\sgn(e(\eta))W_s,
$$
where
$$
W_s=\sum_{k=0}^n\binom nk\binom d{k+5-s}
   =\binom{2d-5}{d-5+s}.
\]
The last equality follows from Vandermonde's identity and symmetry of binomial coefficients. In particular, $W_s=W_{5-s}$. For now we apply this rule also to the two exceptional counts.

Suppose first that every shift with $\eta_1=0$ gave a decrease and every shift with $\eta_1=1$ gave an increase. Pairing each $\eta$ with its coordinatewise complement $\mathbf1-\eta$ would cancel all contributions: the signs are opposite and the weights $W_s$ and $W_{5-s}$ are equal. We must remove the contributions assigned to the five unchanged shifts in $Z$. The two shifts $00110,01001$ were assigned $+W_2$ each, while $10100,11000$ were assigned $-W_2$ each, so removing these four changes nothing. The remaining shift $10000$ was assigned $-W_1$. Removing it leaves precisely
\[
W_1=\binom{2d-5}{d-4}.
\]

The condition $k=0$ determines a unique prefix $\delta=(1,\ldots,1)$, with the empty prefix when $n=0$. Hence each of the two exceptional changes occurs exactly once. Finally, the two exceptional changes $2\to1$ and $1\to2$ merely interchange two entries, so their actual combined contribution is zero. The rule used above assigned them $\binom d1-\binom d0=d-1$. We must therefore subtract $d-1$. We obtain
\[
\Delta(C_d,F_d)=W_1-(d-1)=\binom{2d-5}{d-4}-d+1,
\]
as required.
\end{proof}

\section{Other observation windows}\label{sec:other-windows}

\subsection{Counting identities}\label{sec:counting}

We now allow the observation window to have arbitrary positive integer side lengths. For a fleet $\ka F$, such a window $W$ and $v\in\Z^d$, let $N_{\ka F}(W;v)$ be the number of black cubes contained in $W+v$. We first consider windows whose side lengths are all even or all odd.

\begin{proposition}\label{prop:even-odd-windows}
Let $\ka F$ be a fleet, $v\in\Z^d$ and $k_1,\ldots,k_d$ be positive integers.
\begin{enumerate}[label=\textup{(\alph*)}]
\item If $W=\prod_{i=1}^d[0,2k_i]$, then
\[
N_{\ka F}(W;v)=\prod_{i=1}^d k_i.
\]
Consequently,
\[
\sum_{\varepsilon\in E_d} N_{\mathcal{F}}(W;v+\varepsilon)
=2^d\prod_{i=1}^d k_i
=\mathrm{vol}(W).
\]
\item If $W=\prod_{i=1}^d[0,2k_i-1]$, then
\[
\sum_{\eps\in E_d}N_{\ka F}(W;v+\eps)=\prod_{i=1}^d(2k_i-1)=\mathrm{vol}(W).
\]
\end{enumerate}
\end{proposition}

\begin{proof}
By translating the fleet, we may assume $v=0$. First observe that every cube $[0,2]^d+z$, $z\in\Z^d$, contains exactly one black cube. Indeed, its centre $z+(1,\ldots,1)$ is a vertex of exactly one black cube, since the discrete cubes of the fleet form a partition of $\Z^d$. That black cube is contained in $[0,2]^d+z$. Conversely, every black cube contained in $[0,2]^d+z$ has $z+(1,\ldots,1)$ as a vertex.

\textup{(a)} We divide $W$ into the $\prod_i k_i$ cubes
\[
[0,2]^d+2z,\qquad z\in\prod_{i=1}^d\{0,\ldots,k_i-1\}.
\]
These cubes have pairwise disjoint interiors, and each contains exactly one black cube. Since black cubes have integer vertices and side length $1$, none crosses an integer coordinate hyperplane in its interior. Thus every black cube contained in $W$ belongs to exactly one cube of the partition, giving the asserted count.

\textup{(b)} Now we divide $W$ into the unit cubes
\[
K_z=[0,1]^d+z,\qquad z\in Z:=\prod_{i=1}^d\{0,\ldots,2k_i-2\}.
\]
Let us form a table with rows indexed by $z\in Z$ and columns indexed by $\eps\in E_d$. Put $1$ in position $(z,\eps)$ if $K_z+\eps$ is a black cube, and put $0$ otherwise. For each fixed $z$,
\[
\bigcup_{\eps\in E_d}(K_z+\eps)=[0,2]^d+z.
\]
This cube contains exactly one black cube, which, by integrality of its vertices, equals exactly one of the cubes $K_z+\eps$. Hence every row of the table contains exactly one $1$.

For each fixed $\eps$, the cubes $K_z+\eps$, $z\in Z$, form a partition of $W+\eps$. Every black cube contained in this window equals one of these unit cubes, so the sum of column $\eps$ is $N_{\ka F}(W;\eps)$. Counting the ones by rows and by columns gives
\[
\sum_{\eps\in E_d}N_{\ka F}(W;\eps)=|Z|=\prod_{i=1}^d(2k_i-1).
\]
\end{proof}

When all side lengths are even, every pair of scores is a draw. By part~\textup{(b)}, exactly one of the windows $[0,1]^d+v+\eps$, $\eps\in E_d$, contains a black cube. For $W=[0,3]^d$, the same identity gives another proof of Proposition~\ref{prop:count}.

\begin{corollary}\label{cor:mixed-windows}
Let $\ka F$ be a fleet, $v\in\Z^d$ and $W=\prod_{i=1}^d[0,a_i]$, where $a_i=2k_i$ for $i\in I$ and $a_j=2k_j-1$ for $j\in J=[d]\setminus I$, with all $k_i\ge1$. Then
\[
\sum_{\eps\in E_d}N_{\ka F}(W;v+\eps)=\mathrm{vol}(W).
\]
More precisely, let $E_J=\{\eps\in E_d:\eps_i=0\text{ for }i\in I\}$. Summing only over shifts in the coordinates with odd side lengths, we have
\begin{equation}\label{eq:mixed-window-count}
\sum_{\eps\in E_J}N_{\ka F}(W;v+\eps)=\prod_{i\in I}k_i\prod_{j\in J}(2k_j-1).
\end{equation}
Consequently, the mean of $N_{\ka F}(W;v+\eps)$ over $\eps\in E_d$ equals $\mathrm{vol}(W)/2^d$.
\end{corollary}

\begin{proof}
We first prove~\eqref{eq:mixed-window-count}. As before, we may assume $v=0$. We divide $W$ into boxes $K$ with side length $2$ in the coordinates in $I$ and side length $1$ in the coordinates in $J$, using intervals $[2r,2r+2]$ in the former coordinates and $[s,s+1]$ in the latter. There are $\prod_{i\in I}k_i\prod_{j\in J}(2k_j-1)$ such boxes.

We use the same table argument as in Proposition~\ref{prop:even-odd-windows}\textup{(b)}, now with rows indexed by these boxes and columns by $\eps\in E_J$. The entry in position $(K,\eps)$ is the number of black cubes contained in $K+\eps$. For a fixed $K$, the union of these translates is a cube of side $2$ with integer vertices. Its unique black cube lies in exactly one of the translates $K+\eps$, so each row sums to $1$. For a fixed $\eps$, the translates partition $W+\eps$, and no black cube crosses their boundaries in its interior. The column sum is therefore $N_{\ka F}(W;\eps)$, proving~\eqref{eq:mixed-window-count}.

To obtain the sum over all shifts, we let $E_I$ be the shifts in $E_d$ supported on $I$. We apply~\eqref{eq:mixed-window-count} at each base point $v+\eta$, $\eta\in E_I$, and sum. Each shift in $E_d$ has a unique expression $\eta+\eps$ with $\eta\in E_I$ and $\eps\in E_J$. Thus
\[
\sum_{\eps\in E_d}N_{\ka F}(W;v+\eps)
=2^{|I|}\prod_{i\in I}k_i\prod_{j\in J}(2k_j-1)
=\mathrm{vol}(W).
\]
Dividing by $2^d$ gives the mean.
\end{proof}

If exactly one side length is odd, say $J=\{j\}$, the mixed identity also shows that neither of two fleets beats the other. Indeed, put $\sigma=(2k_j-1)\prod_{i\ne j}k_i$. Two scores whose shifts differ only in coordinate $j$ sum to $\sigma$, so flipping this coordinate pairs each score $x$ with $\sigma-x$, for either fleet. Applying the flip to both players' shifts reverses every strict comparison and preserves every draw. Hence the numbers of wins and losses are equal.

Thus a strict advantage requires at least two odd side lengths. With only one, the scores may still vary, but the symmetry above forces the numbers of wins and losses to agree. As in the study of local number fluctuations~\cite{TorquatoStillinger2003,TorquatoKimKlatt2021}, the mean depends only on the volume of the window, while the distribution of the scores reflects the geometry of the fleet.

\subsection{Changing the size of the window}\label{sec:window-size}

Changing the size of the window can change the winner. We use the periodic realization $\ka F_H$ of a fleet code $H$ from Proposition~\ref{prop:local-code}. For an integer $a\ge1$ and $\eps\in E_d$, let
\[
N_a(H,\eps)=N_{\ka F_H}([0,a]^d;\eps),
\]
and let $L_a(H)$ be the corresponding list of scores. Define $w_a$ and $\Delta_a$ as before, using these scores. In particular, $\Delta_3=\Delta$, while $\Delta_a=0$ for even $a$ by Proposition~\ref{prop:even-odd-windows}.

To compute the scores, write $a=4m+r$, where $0\le r<4$. Every black cube of the periodic fleet has the form $[0,1]^d+h+4z$, with $h\in H$ and $z\in\Z^d$. Such a cube is contained in the observation window precisely when
\[
[0,1]^d+h+4z\subseteq[0,a]^d+\eps
\quad\Longleftrightarrow\quad
\eps_i\le h_i+4z_i\le\eps_i+a-1\quad(i\in[d]).
\]
For a fixed $h\in H$, we therefore count, in each coordinate, the integers congruent to $h_i$ modulo $4$ among $\eps_i,\ldots,\eps_i+a-1$. Each residue modulo $4$ occurs $m$ times, with one additional occurrence for the first $r$ residues in this interval. We use Iverson brackets: $[A]=1$ if the condition $A$ holds, and $[A]=0$ otherwise. Thus the number of choices for $z_i$ is $m+[(h_i-\eps_i)\bmod4<r]$.

The coordinates $z_i$ can be chosen independently. Hence the number of cubes $[0,1]^d+h+4z$ contained in the window is the product of these numbers over all coordinates. Every cube of the periodic fleet has a unique representation of this form, since $H\subseteq\{0,1,2,3\}^d$. Summing over $h\in H$ gives
\begin{equation}\label{eq:radar-count-polynomial}
P_{H,r,\eps}(m):=N_{4m+r}(H,\eps)=\sum_{h\in H}\prod_{i=1}^d\bigl(m+[(h_i-\eps_i)\bmod4<r]\bigr),
\end{equation}
a polynomial in $m$ with leading term $2^dm^d$. 

For example, when $a=5$, we have $N_5(H,\eps)=\sum_{h\in H}2^{|\{i:h_i=\eps_i\}|}$.

\subsection{Possible behavior as the window grows}\label{sec:window-behavior}

\paragraph{Eventual stabilization.}
Fix two fleet codes $F,G$ in dimension $d$ and a residue $r\in\{1,3\}$. For every window size, both players' shifts are chosen from the same fixed set $E_d=\{0,1\}^d$, so there are always $4^d$ pairs of shifts. The balance is
\[
\Delta_{4m+r}(F,G)
=\sum_{(\eps,\delta)\in E_d\times E_d}
\sgn\!\bigl(P_{F,r,\eps}(m)-P_{G,r,\delta}(m)\bigr).
\]
Each summand is $1$ for a win by $F$, $-1$ for a win by $G$, and $0$ for a draw, with $\sgn(0)=0$.

For each fixed pair $(\eps,\delta)$, the difference of the two score polynomials is either identically zero or has a constant strict sign for all sufficiently large $m$. Hence each summand eventually becomes constant. Since there are only $4^d$ summands, they are all constant beyond a common threshold, and so the value of their sum is constant as well. Applying this argument for both values of $r$, we obtain integers $A_1,A_3$ and $m_0$ such that
\begin{equation}\label{eq:radar-eventual-two-values}
\Delta_{4m+1}(F,G)=A_1\quad\text{and}\quad \Delta_{4m+3}(F,G)=A_3\qquad(m\ge m_0).
\end{equation}
Thus stabilization takes place separately along the side lengths $4m+1$ and $4m+3$. The eventual behavior is described by four cases:
\begin{enumerate}[label=\textup{(\alph*)},leftmargin=*]
\item $A_1=A_3=0$: neither code beats the other for any sufficiently large odd side length.
\item $A_1A_3>0$: the same code wins for every sufficiently large odd side length.
\item $A_1A_3<0$: the strict winner alternates between consecutive sufficiently large odd side lengths.
\item Exactly one of $A_1,A_3$ is zero: equal numbers of wins and losses alternate with a strict victory by one code.
\end{enumerate}
All four cases occur in dimension $3$. A zero balance means equal numbers of winning and losing pairs of shifts, and does not require any individual comparison to be a draw.

Each set displayed below consists of eight pairwise dichotomous words and hence is a fleet code; as before, $C_3=\{0,2\}^3$.

The stated balances follow directly from~\eqref{eq:radar-count-polynomial} by comparing the score polynomials for all $64$ pairs of shifts in the indicated ranges of $m$.

\paragraph{Case (a).}
Let
\[
\begin{aligned}
F&=\{000,002,020,022,200,202,221,223\},\\
G&=\{000,002,020,122,201,203,220,322\}.
\end{aligned}
\]
For every $m\ge2$,
\[
\Delta_{4m+1}(F,G)=0,\qquad \Delta_{4m+3}(F,G)=0.
\]
Thus $(A_1,A_3)=(0,0)$.

\paragraph{Case (b).}
Let
\[
H=\{000,002,020,022,211,213,231,233\}.
\]
For every $m\ge1$,
\[
\Delta_{4m+1}(H,C_3)=6,\qquad \Delta_{4m+3}(H,C_3)=6.
\]
Thus $(A_1,A_3)=(6,6)$.

\paragraph{Case (c).}
Let
\[
\begin{aligned}
U&=\{000,002,020,022,201,203,221,223\},\\
V&=\{001,013,021,033,210,212,231,233\}.
\end{aligned}
\]
For every $m\ge2$,
\[
\Delta_{4m+1}(U,V)=4,\qquad \Delta_{4m+3}(U,V)=-4.
\]
Thus $(A_1,A_3)=(4,-4)$.

\paragraph{Case (d).}
Let
\[
K=\{001,003,020,022,201,213,221,233\}.
\]
For every $m\ge1$,
\[
\Delta_{4m+1}(C_3,K)=0,\qquad \Delta_{4m+3}(C_3,K)=-6.
\]
Thus $(A_1,A_3)=(0,-6)$.

Direct enumeration of the two labelled fleet codes in dimension $1$ and the twelve in dimension $2$ gives $(A_1,A_3)=(0,0)$ for every pair. Thus dimension $3$ is the smallest dimension in which any of cases~\textup{(b)--(d)} can occur. The comparison uses the eventual signs of the score-polynomial differences, so it does not depend on testing an arbitrarily chosen large window.

\paragraph{Variation before stabilization in dimension four.}
In dimension $d=4$, let
\[
\begin{aligned}
R=\{&0100,0102,0120,0302,0310,0330,1122,1322,\\
    &2002,2010,2030,2201,2203,2220,3022,3222\},
\end{aligned}
\]
and let
\[
S=\bigl(\{0\}\times\{0,2\}^3\bigr)\cup\bigl(\{2\}\times\{1,3\}^3\bigr).
\]
Both sets consist of sixteen pairwise dichotomous words and hence are fleet codes. Table~\ref{tab:late-stabilization-balances} gives their balances, computed from~\eqref{eq:radar-count-polynomial}; the reverse comparison satisfies $\Delta_a(S,R)=-\Delta_a(R,S)$.
\begin{center}
\small
\begin{tabular}{rrrrr}
\toprule
$m$ & $4m+1$ & $\Delta_{4m+1}(R,S)$ & $4m+3$ & $\Delta_{4m+3}(R,S)$\\
\midrule
0  & 1  & 0  & 3  & $-4$\\
1  & 5  & 13 & 7  & $-12$\\
2  & 9  & 4  & 11 & $-22$\\
3  & 13 & 6  & 15 & $-20$\\
4  & 17 & 8  & 19 & $-12$\\
5  & 21 & 8  & 23 & $-6$\\
6  & 25 & 8  & 27 & $-6$\\
7  & 29 & 0  & 31 & $-4$\\
8  & 33 & 0  & 35 & $-4$\\
9  & 37 & 0  & 39 & $-4$\\
10 & 41 & 0  & 43 & $-4$\\
11 & 45 & 0  & 47 & $-4$\\
12 & 49 & 0  & 51 & $-2$\\
13 & 53 & 0  & 55 & 0\\
\bottomrule
\end{tabular}
\captionof{table}{Balances for the four-dimensional codes $R,S$ before and at stabilization.}
\label{tab:late-stabilization-balances}
\end{center}

The balances vary nonmonotonically even within each residue class modulo $4$. As the odd side length increases from $3$ to $27$, the winner changes at every step, giving twelve strict reversals. From $29$ to $51$, zero balances alternate with victories by $S$. For every odd $a\ge53$, we have $\Delta_a(R,S)=0$, so $(A_1,A_3)=(0,0)$. The eventual identities were verified from~\eqref{eq:radar-count-polynomial} for all $256$ pairs of shifts and all $m\ge13$. Since $\Delta_{51}(R,S)=-2$, the least odd stabilization threshold is exactly $53$. This example motivates the questions below about the size of the stabilization threshold and the number of reversals before it.

\section{Open problems}\label{sec:problems}

\begin{problem}
We conjecture that, for every $d\ge3$, each fleet code in dimension $d$ is beaten by another fleet code in the same dimension. The computations in Section~\ref{sec:low-dimensions} establish the conjecture for $d=3$ and $d=4$.
\end{problem}

\begin{problem}
Determine $\min\{\Delta(F,C_d):F\in\mathscr F_d\text{ non-simple}\}$ for $d\ge5$. Proposition~\ref{prop:unbounded-simple-advantage} gives the upper bound $d-1-\binom{2d-5}{d-4}$, which is negative for every $d\ge5$. More generally, which non-simple codes beat $C_d$? The example $H_5$ shows that allowing changes in three coordinates is already enough for the conclusion of Theorem~\ref{thm:two-coordinate-changes} to fail. What further conditions on the structure of a code ensure that it beats $C_d$?
\end{problem}

\begin{problem}\label{prob:comparability}
Let $p_d$ be the proportion of unordered pairs of distinct labelled codes in $\mathscr F_d$ whose lists are comparable under majorization. We have $p_2=1$, $p_3\approx0.917$ and $p_4\approx0.629$. Does $p_d\to0$ as $d\to\infty$? Which geometric operations on codes, such as suitable twin-pair flips, are monotone with respect to majorization?
\end{problem}

\begin{problem}
Determine, in general, which pairs of eventual values $(A_1,A_3)$ in~\eqref{eq:radar-eventual-two-values} are attainable in each dimension. All four cases described in Section~\ref{sec:window-behavior} occur already in dimension $3$.

For a pair $F,G$, let $a_0(F,G)$ be the least positive odd integer such that $\Delta_a(F,G)=A_r$ whenever $a\ge a_0(F,G)$ is odd and $a\equiv r\pmod4$. How large can $a_0(F,G)$ be in terms of $d$? The four-dimensional example in Section~\ref{sec:window-behavior} has $a_0(R,S)=53$ and exhibits twelve strict reversals before stabilization. Can the number of strict reversals before this threshold be arbitrarily large as $d$ grows? When counting reversals, omit zero balances and count changes between opposite nonzero signs.

For nonperiodic fleets, the polynomial argument is no longer available. For two fixed nonperiodic fleets and a fixed base point, choose both shifts independently and uniformly from $E_d$ for every window size. Can the strict winner change infinitely often while the sequence of outcomes fails to become eventually periodic? Throughout these questions, coordinate and letter labels remain fixed, since the game is not invariant under code isomorphisms.
\end{problem}


\paragraph{Acknowledgments.}
The author thanks Jarek Grytczuk and Krzysiek Przes{\l}awski for helpful discussions.

\section*{Statements and Declarations}

\paragraph{Funding.}
No funding was received for conducting this study.

\paragraph{Competing interests.}
The author has no relevant financial or non-financial interests to disclose.

\paragraph{Data and code availability.}
The code and computational data supporting the results of this study are available from the corresponding author upon reasonable request.

\paragraph{Author contributions.}
Andrzej P. Kisielewicz is the sole author and is responsible for the research and the preparation of the manuscript. The use of AI tools is disclosed in the Introduction.

\end{document}